\documentclass{article}
\usepackage[english, spanish]{babel}
\usepackage[T1]{fontenc}
\usepackage[utf8]{inputenc}
\usepackage{amsmath}
\usepackage{amssymb}
\usepackage{amsthm}
\usepackage{dsfont}
\usepackage{enumitem}
\usepackage{centernot}

\usepackage{etoolbox}
\usepackage{mathtools}
\usepackage{mathrsfs}

\usepackage{hyperref}

\usepackage{calc}

\usepackage{xcolor}
\usepackage{tikz}
\usetikzlibrary{tikzmark,calc}
\usepackage[all,cmtip]{xy}

\newtoggle{useenglish}
\togglefalse{useenglish}

\theoremstyle{plain}
\newtheorem{teo}{\iftoggle{useenglish}{Theorem}{Teorema}}
\newtheorem{lem}[teo]{\iftoggle{useenglish}{Lemma}{Lema}}
\newtheorem{prop}[teo]{\iftoggle{useenglish}{Proposition}{Proposici\'on}}
\newtheorem{cor}[teo]{\iftoggle{useenglish}{Corollary}{Corolario}}

\theoremstyle{definition}
\newtheorem{defi}{\iftoggle{useenglish}{Definition}{Definici\'on}}
\newtheorem{eje}{\iftoggle{useenglish}{Example}{Ejemplo}}
\newtheorem{obs}{\iftoggle{useenglish}{Observation}{Observaci\'on}}
\newtheorem{nota}{\iftoggle{useenglish}{Remark}{Nota}}
\newtheorem*{dem}{\iftoggle{useenglish}{Proof}{Demostración}}

\DeclarePairedDelimiter{\abs}{\lvert}{\rvert}

\newsavebox{\delimpair}

\DeclarePairedDelimiterX{\truthval}[1]{[}{]}{\hspace{-0.4em}\delimsize[{#1}\delimsize]\hspace{-0.4em}}

\newcommand{\printwidth}[1]{\newlength{\dropthis}\settowidth{\dropthis}{#1}\the\dropthis}

\newlength{\strikelengthinternal}

\DeclareMathOperator{\id}{id}

\DeclareMathOperator{\Orb}{Orb}

\DeclareMathOperator{\Sym}{Sym}
\DeclareMathOperator{\Per}{Per}

\DeclareMathOperator{\Aut}{Aut}

\DeclareMathOperator{\Fix}{Fix}

\newcommand{\lang}{\mathscr{L}}
\newcommand{\Alf}{\mathcal{A}}

\newcommand{\R}{\mathds{R}}
\newcommand{\N}{\mathds{N}}
\newcommand{\Z}{\mathds{Z}}

\newcommand{\GL}{\mathrm{GL}}

\newcommand{\shift}[1]{\mathsf{X}_{#1}}

\newcommand{\actson}[1][ ]{\mathrel{\stackrel{#1}{\curvearrowright}}}

\newcommand{\dfn}{\mathrel{\raisebox{0.06ex}{$:$}{=}}}

\newcommand{\ind}{\mathds{1}}

\newcommand{\dotcup}{\mathbin{\ensuremath{\mathaccent"7201\cup}}}

\newcommand{\incl}{\mathrel{{\hookrightarrow}}}
\newcommand{\epi}{\mathrel{{\twoheadrightarrow}}}

\newsavebox{\fmbox}

\toggletrue{useenglish}

\begin{document}
	\selectlanguage{english}
	
	%\title{Computation of extended symmetry groups for multidimensional subshifts with hierarchical structure}
	\title{Extended symmetry groups of multidimensional subshifts with hierarchical structure}
	\author{Álvaro Bustos\footnote{Contact e-mail: \texttt{abustos@dim.uchile.cl}. The author was supported by CONICYT Doctoral Fellowship 21171061(2017).}}
	\date{Departamento de Ingeniería Matemática \\ Universidad de Chile \\ Beauchef 851, Santiago, Chile \\\vspace*{1em} January 23, 2019}
		
	\maketitle
	
	\begin{abstract}

		We study the automorphism group, i.e. the centralizer of the shift action inside the group of self-homeomorphisms, together with the extended symmetry group (the corresponding normalizer) of certain $\Z^d$ subshifts with a hierarchical structure like bijective substitutive subshifts and the Robinson tiling.
		Treating those subshifts as geometrical objects, we introduce techniques to identify allowed symmetries from large-scale structures present in certain special points of the subshift, leading to strong restrictions on the group of extended symmetries. We prove that in the aforementioned cases, $\Sym(X, \Z^d)$ (and thus $\Aut(X, \Z^d)$) is virtually-$\Z^d$ and we explicitly represent the nontrivial extended symmetries, associated with the quotient $\Sym(X, \Z^d)/\Aut(X, \Z^d)$, as a subset of rigid transformations of the coordinate axes. We also show how our techniques carry over to the study of the Robinson tiling, both in its minimal and non minimal version.

	\end{abstract}
	
	\section{Preliminaries}
	
	This section will introduce the basic concepts and notation to be used in what follows. We assume some basic familiarity with symbolic dynamics (namely, the nature of a shift space and the shift action); for the reader interested in more in-depth treatment of this subject, we recommend consulting the book by Lind and Marcus \cite{LM95} for the one-dimensional case and the text by Ceccherini-Silberstein and Coornaert \cite{TCS2010} as an introduction to the case of general groups. 
	
	The subject of symbolic dynamics deals with a specific kind of group action, the \textbf{shift action:}
	
	\begin{defi}
		Let $\mathcal{A}$ be a finite set (which we shall call \textbf{alphabet}) and $G$ any group (usually, but not always, assumed to be finitely generated). The \textbf{full-shift} is the topological space $\mathcal{A}^G$, with the product topology (after giving $\mathcal{A}$ the discrete topology). The (left) \textbf{shift action} is the following group action $G\actson[\sigma]\mathcal{A}^G$:
		\[(\forall x=(x_g)_{g\in G}\in\mathcal{A}^G)(\forall g,h\in G):(\sigma_g(x))_h\dfn x_{g^{-1}h}. \]
		A closed subset of $\mathcal{A}^G$ that is invariant under the shift action is called a $G$-\textbf{subshift}.
	\end{defi}
	
	Subshifts are usually defined combinatorially instead of topologically, via \textbf{forbidden patterns}. A \textbf{pattern} $P$ with (finite) \textbf{support} $U\subset G$ is a function $P:U\to\mathcal{A}$. $G$ acts on the set of all patterns $\mathcal{A}^{*,G}$ by translation: the pattern $g\cdot P$ has support $gU$ and is defined by $(g\cdot P)_{h}=P_{g^{-1}h}$. We say that a point $x\in\mathcal{A}^G$ \textbf{contains} the pattern $P$ (usually written as $P\sqsubset x$) if, for some $g\in G$, $x|_{gU}=g\cdot P$; note that, by this definition, any translation of $P$ is contained in $x$ as well, and thus is functionally ``the same'' as $P$. A similar definition applies for two patterns $P$ and $Q$, where we again use $P\sqsubset Q$ to denote the subpattern relation.
	
	Any subshift $X$ can be described by a set of forbidden patterns $\mathcal{F}\subseteq\mathcal{A}^{*,G}$; namely, given such a set $\mathcal{F}$ we define the following subset of $\mathcal{A}^G$:
	\[\shift{\mathcal{F}}\dfn\{x\in\mathcal{A}^G:(\forall P\in\mathcal{F}):P\not\sqsubset x \}. \]
	It is not hard to prove that $\shift{\mathcal{F}}$ is a subshift and that any subshift $X$ is equal to $\shift{\mathcal{F}}$ for some (usually not uniquely determined) set of patterns $\mathcal{F}$. If the set $\mathcal{F}$ can be chosen finite, we say that $X$ is a \textbf{shift of finite type} (or $G$-SFT).
	
	Since subshifts can be regarded as both dynamical and combinatorial objects, we can classify them not only in combinatorial terms (shifts of finite type, substitutive subshifts, etc.) but also in terms of their dynamics. One of the main classifications we shall be interested in is given by the following definition:
	\begin{defi}
		Let $X\subseteq\mathcal{A}^{G}$ be a $G$-subshift. We say that the shift action $G\actson[\sigma]X$ is \textbf{faithful} if for all $g\in G\setminus\{1_G\}$ there is a point $x\in X$ such that $\sigma_g(x)\ne x$, i.e. if $\sigma_g=\id$ implies $g=1_G$.
	\end{defi}
	We shall assume unless otherwise stated that the shift action is faithful for the shifts under study. This is because, if the shift action is not faithful, there is a strict subgroup $H<G$ such that $X$ essentially behaves like an $H$-subshift; thus, we can always limit ourselves to the faithful case. Moreso, in the study of the automorphism group described below, having a faithful action makes it easier to describe such a group.
	
	Let $X,Y$ be $G$-subshifts. A continuous, shift-commuting ($f\circ\sigma_g=\sigma_g\circ f$ for all $g\in G$) map $f:X\to Y$ is called a \textbf{sliding block code}. These kinds of mappings act as structure-preserving morphisms for this class of group actions (or dynamical systems). In particular, when there is a bijective sliding block code from $X$ to $Y$, both subshifts share all topological and dynamical properties such as periodic points, isolated points, dense subsets, etc.; thus, we call such a mapping an \textbf{isomorphism} (or, when $X=Y$, \textbf{automorphism}) and $X$ and $Y$ \textbf{isomorphic} subshifts. The name ``sliding block code'' comes from the following well-known result:
	
	\begin{teo}[Curtis-Hedlund-Lyndon]
		Let $X$ and $Y$ be $G$-subshifts over finite alphabets $\mathcal{A}_X$ and $\mathcal{A}_Y$, respectively. A mapping $\varphi:X\to Y$ is a sliding block code if, and only if, there exists a finite subset $\mathcal{U}\subset G$ (called the \textbf{window} of $\varphi$) and a mapping $\Phi:\mathcal{A}_X^U\to\mathcal{A}_Y$ (called the \textbf{local function} associated to $\varphi$) such that:
		\[\varphi(x)_g = \Phi(x|_{gU}), \]
		in which we identify the pattern $x|_{gU}$ with its corresponding translate $g^{-1}\cdot x|_{gU}$ with support $U$.
	\end{teo}
	
	\begin{obs}
		Let $\varphi:X\to Y$ be a sliding block code given by a local function $\Phi:\mathcal{A}_X^U\to\mathcal{A}_Y$. If $V\supset U$ is a larger finite set, we may define a new local function $\Phi':\mathcal{A}_X^V\to\mathcal{A}_Y$ which induces the same sliding block code $\varphi$ (this can be seen by taking $\Phi'(P)=\Phi(P|_U)$ for all $P\in\mathcal{A}_X^V$). Thus, in $\Z^d$, we can always assume that the set $U$ is of the form $[-r\vec{\ind},r\vec{\ind}]$ for some $r\ge 0$; the least $r$ for which there is a local function $\Phi$ with a window of this form for $\varphi$ is often called the \textbf{radius} of the sliding block code $\varphi$. A sliding block code of radius $0$ is often called a \textbf{relabeling map}.
	\end{obs}
	
	Our main subject of study, at least for the first part of this work, is the set of automorphisms $f:X\to X$, which we shall denote as $\Aut(X,G)$; note that this set is a group under the operation of composition, and that $Z(G)$, the center of the group $G$, embeds into $\Aut(X,G)$ because, if $gh=hg$, then $\sigma_g\circ\sigma_h=\sigma_{gh}=\sigma_{hg}=\sigma_h\circ\sigma_g$ and thus for any $h\in Z(G)$ the mapping $\sigma_h$ is a continuous, shift-commuting homeomorphism. This is very important in the case of abelian groups such as $\Z^d$, where $G=Z(G)$.
	
	In what follows, we shall be mostly concerned with free abelian groups, namely $\Z^d$; in the following text, unless stated otherwise, we reserve the letter $d$ for the rank or dimension of this underlying group. We will denote elements of this group with vector notation, $\vec{k}=(k_1,\dots,k_d),k_1,\dots,k_d\in\Z$. In this context, the letter $\vec{s}$ will be used for a specific, fixed ``size'' number (as we shall see below) and we will use other letters such as $\vec{k}$ and $\vec{p}$ for generic elements of $\Z^d$.
	
	We shall be mostly concerned with substitutive subshifts, at least in the early sections of this work. These subshifts come from a partial generalization of the concept of a one-dimensional substitution, which consists of a function $\theta:\mathcal{A}\to\mathcal{A}^*\setminus\{\varnothing\}$ that assigns a (nonempty) word comprised of symbols from the alphabet $\mathcal{A}$ to every symbol $a\in\mathcal{A}$. This function extends to the whole of $\mathcal{A}^*$ by concatenation, i.e. we define $\theta^*:\mathcal{A}^*\to\mathcal{A}^*$ by:
	\[\theta^*(a_1a_2\dotsc a_k)=\theta(a_1)\theta(a_2)\dotsc\theta(a_k), \]
	and in a similar fashion we extend $\theta$ to infinite and bi-infinite sequences from $\mathcal{A}^\N$ and $\mathcal{A}^\Z$, respectively. We shall assume the substitution function $\theta$ to be \textbf{primitive}, i.e. there is a fixed $k\in\N$ such that $(\forall a\in\mathcal{A}):\theta^k(a)$ contains all of the symbols of the alphabet $\mathcal{A}$. By taking any (fixed) $x\in\mathcal{A}^\Z$ and applying the substitution repeatedly, we obtain a sequence of points $x^{(n+1)}=\theta_\infty(x^{(n)})$ which, under the mild condition of a primitive substitution being primitive, converges either to a point $x$ that is fixed under $\theta_\infty$ or to a finite, periodic orbit under $\theta_\infty$. Taking the orbit closure (under the shift action) of such a fixed or periodic point, we define a subshift with interesting properties, which is called the \textbf{substitutive subshift} associated to $\theta$.
	
	In the multidimensional case, we may consider a substitution as a mapping that assigns a pattern from $\mathcal{A}^{*,G}$, to each symbol $a\in\mathcal{A}$; however, for arbitrary patterns it is hard to describe the extension of this mapping to any configuration, as there is no direct analogue to concatenation. Thus, we restrict ourselves to the case in which all patterns $\theta(a),a\in\mathcal{A}$ share the same rectangular support $S=[\vec{0},\vec{s}-\vec{\ind}]\dfn\prod_{k=1}^{d}[0,s_k-1]$, and thus $\theta$ is a mapping $\mathcal{A}\to\mathcal{A}^S$. This is called a \textbf{rectangular substitution}. This kind of substitution has the obvious advantage of the concatenation rule being easy to describe: symbols $a_1$ and $a_2$ adjacent in the $\vec{e}_j$ direction result in patterns $\theta(a_1)$ and $\theta(a_2)$ appearing adjacent in the same direction $\vec{e}_j$.
	
	\begin{figure}[ht]
		\centering
		\includegraphics[scale=0.8]{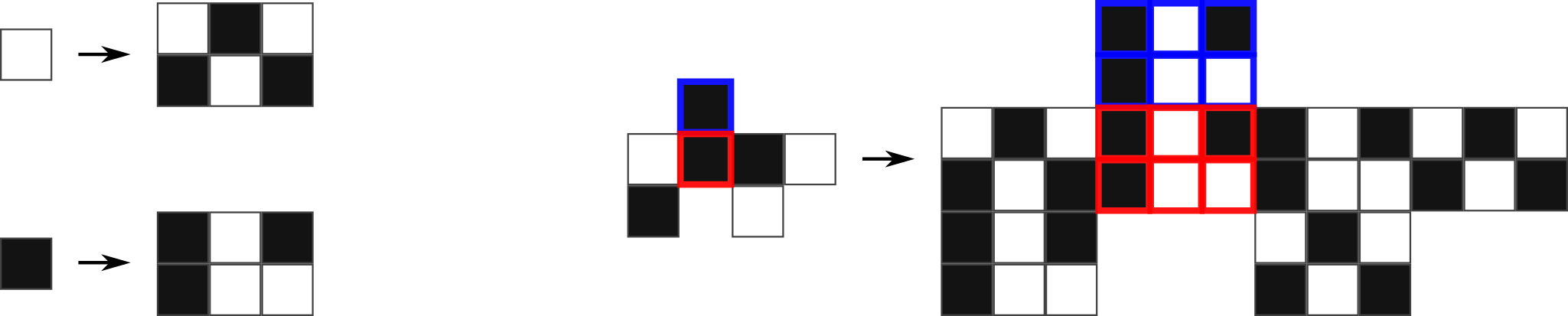}
		\caption{An example of application of a rectangular substitution to a pattern.}
	\end{figure}
	Formally, the extension of $\theta:\mathcal{A}\to\mathcal{A}^S$ to a function $\theta_\infty:\mathcal{A}^{\Z^d}\to\mathcal{A}^{\Z^d}$ follows the same principle: in $\theta_\infty(x)$, every symbol $x_{\vec{k}}$ is replaced by a pattern $\theta(x_{\vec{k}})$, keeping adjacencies, and thus:
	\[\theta_\infty(x)_{\vec{m}\cdot\vec{s}+\vec{k}}=\theta(x_{\vec{m}})_{\vec{k}},\quad\vec{m}\in\Z^d,\vec{k}\in S, \]
	
	Here and in what follows, the multiplication $\vec{m}\cdot\vec{s}$ is taken to be componentwise. For simplicity, we shall always assume that the vector $\vec{s}$ determining the shape of the set $S$ satisfies the condition $\min(s_1,\dots,s_d)>1$; otherwise the problem reduces to the study of a lower-dimensional substitution, and in the one-dimensional case the failure of this condition makes the substitution trivial (in particular, it cannot be primitive). The condition $s_i>1$ for all $1\le i\le d$ over the size vector $\vec{s}$ of the support $S$ ensures that every $\vec{m}\in\N^d$ is in the support of some $\theta^k$ for sufficiently large $k$.

	We shall only use the subscript $\infty$ in specific situations that might lead to confusion, otherwise distinguishing both functions by context. Similarly, we may use the symbol $\theta^*$ for the extension of $\theta$ to the set of all patterns $\mathcal{A}^{*,G}$, but only when needed.
	
	The above definition allows us to introduce the type of subshifts we intend to study:
	
	\begin{defi}
		Let $\theta$ be a primitive rectangular substitution on the alphabet $\mathcal{A}$, and take $\Sigma$ to be the limit set of $\mathcal{A}^{\Z^d}$ under $\theta_\infty$, i.e. the set of all accumulation points of the sequences $(\theta_\infty^k(x))_{k\in\N}$ for all $x\in \mathcal{A}^{\Z^d}$. Note that this set is actually finite, as the accumulation points of such a sequence depend only on the finite subpattern $x|_{\{-1,0\}^d}$ and each sequence has a finite number of such points. We define the \textbf{substitutive subshift} associated to $\theta$ as the shift-orbit closure of $\Sigma$, that is:
		\[\shift{\theta}\dfn\bigcup_{x\in\Sigma}\overline{\Orb_\sigma(x)}; \]
		additionally, we define the \textbf{minimal substitutive subshift} $\shift{\theta}^\circ$ as the following subset of $\mathcal{A}^{\Z^d}$:
		\[\shift{\theta}^\circ\dfn\{x\in\mathcal{A}^{\Z^d}: (\forall U\subset\Z^d,|U|<\infty)(\exists k\in\N)(\exists a\in\mathcal{A}):x|_U\sqsubset\theta^k(a)  \}. \]
	\end{defi}
	
	\begin{nota}
		Note that the usual definition of substitutive subshift in several sources corresponds to $\shift{\theta}^\circ$, and thus it is often assumed that a substitutive subshift is minimal. However, in the analysis below we need the slightly expanded definition given above, which consists of all the points obtained from a ``seed'' (a pattern with support $\{-1,0\}^d$), their shifts and the corresponding limit points; it is easy to see that $\shift{\theta}^\circ\subseteq\shift{\theta}$, with equality only when every possible seed appears as a subpattern of $\theta^k(a)$ for some $a$.
	\end{nota}
	
	\begin{obs}
		Due to the primitivity, the condition $(\exists k\in\N):x|_U\sqsubset\theta^k(a)$ does not depend on the chosen $a$; moreso, this means that any pattern $x|_U$ from a point $x\in\shift{\theta}^\circ$ appears in some $\theta^m(a)$ for any $a\in\mathcal{A}$ and some sufficiently large $m$. As we shall see below, all points from $\shift{\theta}$ (and thus, from $\shift{\theta}^\circ$) are of the form $\sigma_{\vec{k}_m}(\theta_\infty^m(x))$ for some $x\in\shift{\theta}$ and all values of $m$, and thus are concatenations of patterns of the form $\theta^m(a),a\in\mathcal{A}$; this implies that any pattern of the form $x|_U$ with $x\in\shift{\theta}^\circ$ appears as a subpattern of any other $y\in\shift{\theta}$, and thus $\Orb_\sigma(y)$ is dense in $\shift{\theta}^\circ$ for any $y\in\shift{\theta}^\circ$, i.e. this subshift is minimal in a dynamical sense, justifying the name given above. Since this subshift is a subset of $\shift{\theta}$, the latter is only minimal when it equals $\shift{\theta}^\circ$, since dynamical minimality implies minimality by inclusion among closed, shift-invariant subsets.
	\end{obs}
	
	\begin{eje}
		An example of a subshift arising from a rectangular substitution is the two-dimensional Thue-Morse substitutive subshift, given by the following $\theta_{\rm TM}:\mathcal{A}\to\mathcal{A}^{\{1,2\}^2}$:
		\begin{center}
			\begin{tikzpicture}[scale=0.6]
			\draw (0,-0.5) rectangle (1,0.5) (2,-1) rectangle (4,1);
			\draw[xshift=6cm] (2,-1) rectangle (4,1);
			\fill (2,0) rectangle (3,1) (3,0) rectangle (4,-1);
			\fill[xshift=6cm] (0,-0.5) rectangle (1,0.5) (2,0) rectangle (3,-1) (3,0) rectangle (4,1);
			\node at (1.5,0) {$\mapsto$};
			\node at (7.5,0) {$\mapsto$};
			\end{tikzpicture}
		\end{center}
		in which the alphabet $\mathcal{A}$ consists of black and white tiles (identified with the symbols $1$ and $0$, respectively). The $d$-dimensional analogue $\theta_{\rm TM}$ sends each $i\in\{0,1\}$ to a $2\times\cdots\times 2$ pattern with support $S=\{0,1\}^d$ given by:
		\[\theta_{\rm TM}(i)_{(k_1,\dots,k_d)}=\begin{cases}
		0 & \text{if }i+k_1+\dots+k_d\equiv 0\pmod{2}, \\
		1 & \text{otherwise}.
		\end{cases} \]
		In the one-dimensional case, $\shift{\theta_{\rm TM}} = \shift{\theta_{\rm TM}}^\circ$ and a typical point of this subshift looks like the following:
		\[\dotsc 1001011001101001{.}0110100110010110  \dotsc \]
		
		In Figure \ref{fig:contradiction_points} (see page \pageref{fig:contradiction_points})  we can see fragments from two points of the shift $\shift{\theta_{\rm TM}}$, in the case $d=2$. The left one belongs to the minimal subshift $\shift{\theta_{\rm TM}}^\circ$, while the point corresponding to the figure on the right does not.
	\end{eje}
	
	As stated above, we shall refer to any pattern given by a mapping $\{-1,0\}^d\to\mathcal{A}$ as a \textbf{seed}. Any periodic (w.l.o.g. fixed, by replacing $\theta$ by a suitable power $\theta^m$) point of the substitution $\theta$ is uniquely determined (at least in the primitive case with nontrivial support) by its finite subconfiguration with support $\{-1,0\}^d$ and thus by a unique seed.
	
	To work with iterated substitutions and automorphisms, the following notations are useful:
	\begin{align*}
		S^{(m)} &\dfn [\vec{0},\vec{s}^m-\vec{\ind}] = \{\vec{k}=(k_1,\dots,k_d)\in\Z^d: 0\le k_j< s_j^m, 1\le j\le d \},\\
		R^{\circ m} &\dfn \{\vec{k}\in R: [\vec{k}-m\vec{\ind},\vec{k}+m\vec{\ind}]\subseteq R  \}, \text{ for any }R\subseteq\Z^d.
	\end{align*}
	$S^{(m)}$ is the iterated componentwise multiplication of the elements of the set $S$ with themselves, repeated $m$ times. When $\theta$ is a rectangular substitution with support $S$, the set $S^{(m)}$ is the support of any pattern $\theta^m(a),a\in\mathcal{A}$; we may take $\theta^m=(\theta^*)^m|_\mathcal{A}$ as a new substitution that induces the same substitutive subshift as $\theta$ (in the primitive case), which sometimes proves useful to simplify arguments. In the same way, $R^{\circ m}$ refers to the subset of $R$ comprised of all elements ``at distance at least $m$'' from the complement of $R$, i.e. a sort of ``interior''  of $R$. It is easy to see that, if $f$ is a sliding block code of radius $m$, then for any subset $R\subseteq\Z^d$ we have that $x|_R=y|_R$ implies $f(x)|_{R^{\circ m}}=f(y)|_{R^{\circ m}}$.
	
	In the following text we have to deal with the group of $p$-adic numbers, defined as an inverse limit of a chain of cyclic groups, as follows:
	
	\begin{defi}
		Let $p>1$ be a fixed positive integer (not necessarily prime). The \textbf{group of $p$-adic integers} $\Z_p$ is the inverse limit associated to the following diagram of groups and group morphisms:
		\[\xymatrix{\Z/p\Z & \Z/p^2\Z\ar[l]_-{\pi_1} & \Z/p^3\Z\ar[l]_-{\pi_2} & \dots\ar[l]_-{\pi_3}}\]
		where each $\pi_i:\Z/p^{i+1}\Z\to\Z/p^i\Z$ is the remainder modulo $p^i$ function, i.e. $\pi_i([ap^i+b]_{p^{i+1}})=[b]_{p^i}$. Alternatively, $\Z_p$ corresponds to the following subgroup of the infinite product $\prod_{k=1}^\infty \Z/p^k\Z$:
		\[\Z_p\dfn\left\{([m_k]_{p^k})_{k\ge 1}\in\prod_{k=1}^\infty\Z/p^k\Z : (\forall k\ge 1): m_k\equiv m_{k+1}\pmod{p^k}  \right\}. \]
	\end{defi}
	As a subgroup of an infinite product, addition in $\Z_p$ is performed componentwise. It is easy to see that the sequence $([1]_{p^k})_{k\ge 1}$ belongs to $\Z_p$ and generates an infinite cyclic group, which we identify with $\Z$ (and thus we identify the sequence $([1]_{p^k})_{k\ge 1}$ with the integer $1$). The set $\Z_p\setminus\Z$ is nonempty. For instance, the sequence $([1]_2,[1]_4,[5]_8,[5]_{16},[21]_{32},[21]_{64},\dots)$ belongs to $\Z_2$ but it does not represent an integer.
	
	$\Z_p$ has a topological group structure, induced by the prodiscrete topology on the space $\prod_{k=1}^{\infty}\Z/p^k\Z$; we can define a metric which is analogous to the shift metric in $\mathcal{A}^\N$. With this structure, the example given of a sequence in $\Z_2\setminus\Z$ corresponds to the infinite sum $\sum_{k=1}^{\infty}2^{2k-1}$.
	
	\begin{defi}
		The \textbf{$p$-adic odometer} is the topological dynamical system $(\Z_p,\omega_p)$, where $\omega_p(x)=x+1$ (here, as above, $1$ represents the infinite sequence $([1]_{p^k})_{k\ge 1}$) and $\Z_p$ is taken with the prodiscrete topology.
	\end{defi}
	
	As we shall see, under certain hypotheses the maximal equicontinuous factor of a $d$-dimensional substitutive subshift is a product of odometers. Thus, we introduce the notation $\Z_{(s_1,\dots,s_d)}$ (or simply $\Z_{\vec{s}}$) for the product $\Z_{s_1}\times\dotsm\times\Z_{s_d}$, and identify $\Z^d$ with the corresponding subgroup of $\Z_{\vec{s}}$.
	
	\section{Substitution encodings}
	
	For details on the propositions referenced in this section, the survey by Frank \cite{Fr2005} may be consulted. In the one-dimensional case, the book by K\r{u}rka \cite{Kur2003} gives a treatment of this encoding factor as well. More details about substitutions can be found in the book by Pytheas Fogg \cite{Fog2002} (for the one-dimensional case); the book by Michael Baake and Uwe Grimm \cite{BG2013} gives a good treatment of the multidimensional case as well.
	
	In what follows, let $\theta:\mathcal{A}\to\mathcal{A}^S$ be a rectangular primitive substitution of constant size $\vec{s} = (s_1,\dots,s_d)$, with $S=[\vec{0},\vec{s}-\vec{\ind}]$. We assume that $s_i>1$ for all $1\le i\le d$.
	\begin{lem}
		\label{lem:codified_system_subst}
		Given $x\in\shift{\theta}$, there is a unique $\vec{k}_1\in S$ and $y\in\shift{\theta}$ such that $x=\sigma_{\vec{k}_1}(\theta(y))$. By iterating this process, there is a unique $\vec{k}_m\in S^{(m)}$ such that $x=\sigma_{\vec{k}_m}(\theta^m(y))$ and $\vec{k}_m\equiv\vec{k}_r\pmod{\vec{s}^r}$ for $m>r$, with module equivalence taken componentwise.
	\end{lem}
	This allows us to assign a sequence $([\vec{k}_m]_{\vec{s}^m})_{m\ge 1}$ in $\Z_{\vec{s}}$ to each element of $\shift{\theta}$, and it is easy to see that the sequence associated to $\sigma_{\vec{k}}(x)$ is the sequence assigned to $x$, plus $\vec{k}$. This observation leads to the following known result:
	\begin{prop}
		\label{prop:MEF_of_a_subst_shift}
		The maximal equicontinuous factor of a nontrivial substitutive subshift $\shift{\theta}$ given by a $d$-dimensional substitution of constant length $\theta:\mathcal{A}\to\mathcal{A}^S$ is the product system of $d$ odometers, $(\Z_{s_1},\omega_{s_1})\times\dotsm\times(\Z_{s_d},\omega_{s_d})$ (or equivalently, the group $\Z_{\vec{s}}$ with its subgroup $\Z^d$ acting by addition). The factor morphism $\varphi:\shift{\theta}\to\Z_{\vec{s}}$ sends each $x\in\shift{\theta}$ to the uniquely determined sequence belonging to $\Z_{\vec{s}}$ by the previous lemma.
	\end{prop}
	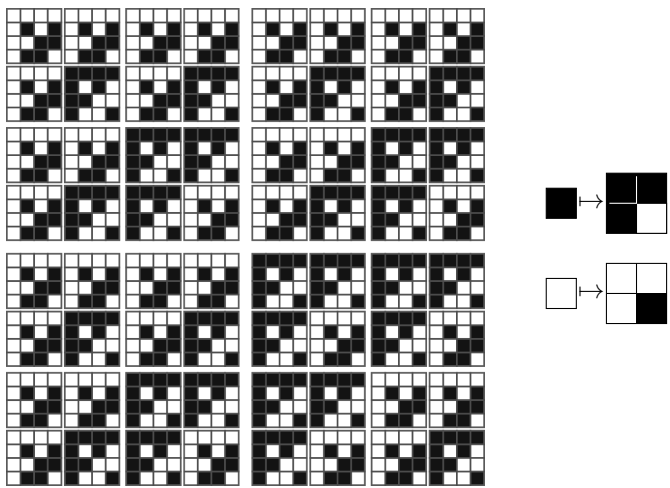
\begin{figure}[ht]
		\centering
		\includegraphics[scale=0.4]{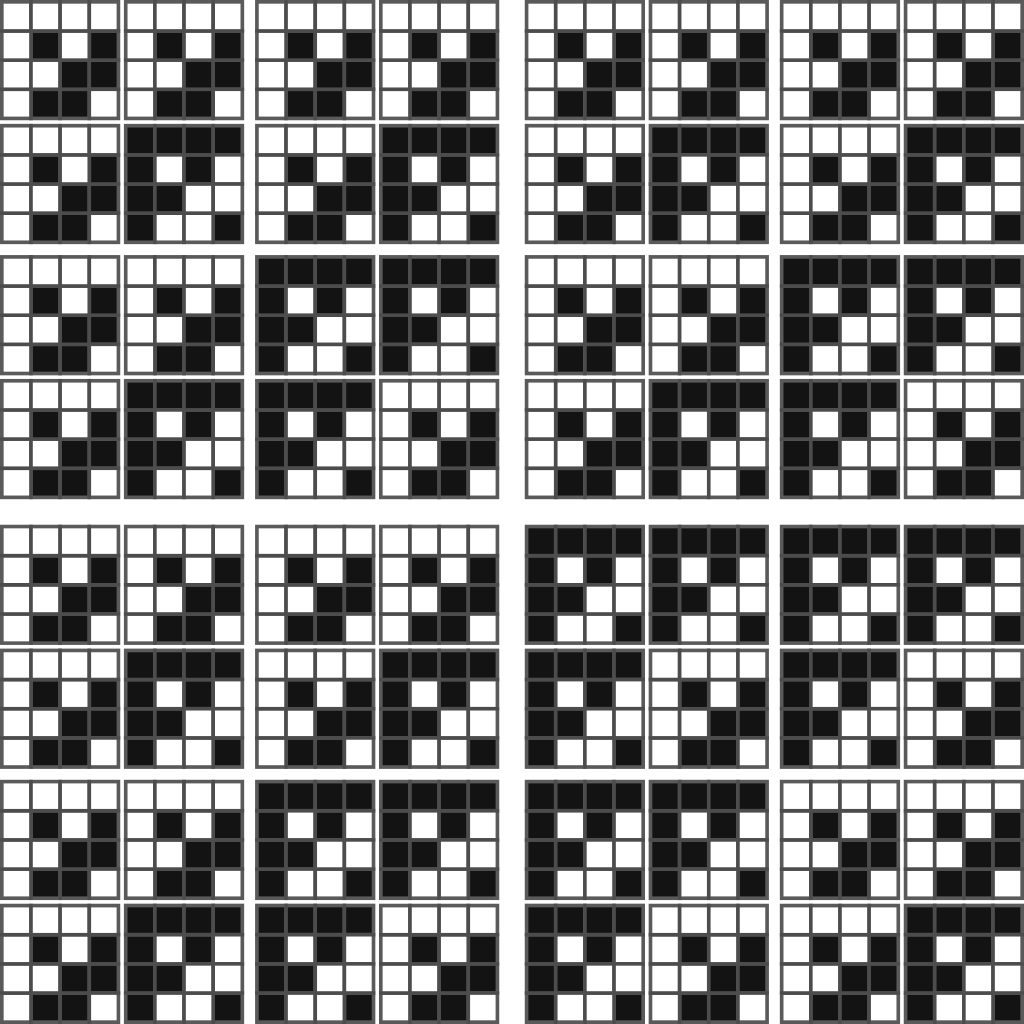}
		\qquad
		\raisebox{2.15cm}{\begin{tikzpicture}[scale=0.4]
			\draw (0,-0.5) rectangle (1,0.5) (2,-1) rectangle (4,1);
			\fill (3,0) rectangle (4,-1);
			\draw (2,0) -- (3,0) -- (3,1);
			\node at (1.5,0) {$\mapsto$};
			\begin{scope}[yshift=3cm]
				\filldraw[fill=black,draw=black] (0,-0.5) rectangle (1,0.5) (2,-1) rectangle (4,1);
				\filldraw[fill=white,draw=black] (3,0) rectangle (4,-1);
				\draw[color=white] (2.1,0) -- (3,0) -- (3,0.9);
				\node at (1.5,0) {$\mapsto$};
			\end{scope}
		\end{tikzpicture}}
		\caption{$2^n\times 2^n$ grids associated with the iterates of a substitution $\theta$ in a point from a substitutive subshift. The corresponding substitution is indicated in the figure.}
	\end{figure}
	
	Proposition \ref{prop:MEF_of_a_subst_shift} is proved in the survey by Frank, \cite{Fr2005}. Note that, in particular, $\varphi$ is continuous; intuitively, if two points $x$ and $y$ match on $[-r\vec{\ind},r\vec{\ind}]$ for sufficiently large $r$, this central pattern is enough to determine the shift from Lemma \ref{lem:codified_system_subst} for all ``small'' values of $m$, and thus $\varphi(x)$ and $\varphi(y)$ match on their first few entries, and thus are ``close'' for the $p$-adic metric.
	
	From the previous proposition, we can prove the following result:
	\begin{lem}
		\label{lem:integer_images_in_MEF}
		Let $x\in\shift{\theta}$ be such that $\varphi(x)=\vec{0}$. Then $x$ is a periodic (w.l.o.g. fixed) point of the substitution $\theta$; the converse is also true. Thus (w.l.o.g. assuming all periodic points of $\theta$ to be fixed points):
		\[\varphi^{-1}[\Z^d]=\bigcup_{x\in\Fix(\theta)}\Orb_{\sigma}(x) \]
	\end{lem}
	
	\begin{dem}
		Without loss of generality we may assume that all periodic points of $\theta$ are fixed points, replacing $\theta$ by a suitable power $\theta^m$ if necessary.
		
		If $\varphi(x)=\vec{0}=([0]_{\vec{s}},[0]_{\vec{s}^2},\dots)$, this means there exists a sequence of points $x_1,x_2,x_3,\dots\in\shift{\theta}$ such that $x=\theta^k(x_k)$, for all $k\ge 1$. Let $P_k$ be the seed obtained from $x_k$ by restriction to the set $\{-1,0\}^d$, for each $k$. Since there are a finite number of possible seeds, there must be a seed $P$ such that $P=P_k$ for infinitely many values of $k$. Notice also that, for any given seed $P$ and any point $y\in\shift{\theta}$ such that its restriction to $\{-1,0\}^d$ is $P$ the sequence $\theta^k(y)$ converges to the same fixed point of $\theta$, $z_P$.
		
		Thus, since $x=\theta^k(x_k)$ for all those infinitely many values of $k$ that have $x_k|_{\{0,1\}^d}=P$, we have that the restriction of $x$ to $[-\vec{s}^k,\vec{s}^k-\vec{\ind}]$ equals the corresponding restriction of $\theta^k(y)$, where $y$ is any point with seed $P$, and thus the distance between $x=\theta^k(x_k)$ and $\theta^k(y)$ is less than $2^{-\min(\vec{s}^k)}$, which goes to zero as $k$ goes to infinity. Thus, by a triangular inequality argument, $x$ must equal the limit $\theta^k(y)$ and thus be a fixed point of $\theta$. Determining the form of $\varphi^{-1}[\Z^d]$ is then direct from the fact that $\Orb_{(\omega_{s_1},\dots,\omega_{s_d})}(\vec{m})=\vec{m}+\Z^d$ in a product of odometers.\qed
	\end{dem}
	
	The importance of these previous results is twofold: the equicontinuous factor $\varphi:\shift{\theta}\epi\Z_{\vec{s}}$ gives a description of a point $x\in\shift{\theta}$ as a concatenation of rectangular patterns of the form $\theta^k(a)$ in a specific way for all values of $k$, and allows us to distinguish certain points, e.g. points with ``fractures''.
	
	\section{Bijective substitutions and Coven's theorem}
	
	We shall be interested in a certain type of rectangular substitutions over the two-symbol alphabet $\mathcal{A}=\{0,1\}$:
	
	\begin{defi}
		A substitution $\theta:\mathcal{A}\to\mathcal{A}^S$ is called \textbf{bijective} if, for each $\vec{k}\in S$ and any $a\ne b\in\mathcal{A}$, we have $\theta(a)_{\vec{k}}\ne\theta(b)_{\vec{k}}$, i.e. each function $\theta_{\vec{k}}:a\mapsto\theta(a)_{\vec{k}}$ is a bijection $\mathcal{A}\to\mathcal{A}$, for all $\vec{k}\in S$.
	\end{defi}
	
	Given a pattern $P$ over the alphabet $\{0,1\}$ with support $L\subset\Z^d$, we write $\overline{P}$ for the pattern with the same support obtained by replacing all $1$s by $0$s and vice versa. It is easy to see that a substitution over the alphabet $\{0,1\}$ is bijective if and only if $\theta(1)=\overline{\theta(0)}$; thus, for any pattern $P$, $\theta(\overline{P})=\overline{\theta(P)}$, and in particular $\theta^k(1)=\overline{\theta^k(0)}$ by induction. Similarly, if we define the substitution $\overline{\theta}$ by the relation $\overline{\theta}(a)=\overline{\theta(a)}$, we can see that $\overline{\theta}^2(a)=\theta^2(a)$ (and hence $\shift{\theta}=\shift{\overline{\theta}}$) so that we may always assume $\theta(a)_{\vec{\ind}}=a,a\in\{0,1\}$. This also implies that for any seed $P$ there exists a periodic point $x_P$ of $\theta$ with seed $P$, and thus $\theta$ has exactly $2^{2^d}$ periodic points (w.l.o.g. we may assume these periodic points to be fixed points by replacing $\theta$ with a suitable power $\theta^m$ such that the symbol on each of the $2^{d}$ corners of the pattern $\theta^m(a)$ is $a$ itself).
	
	Our goal is to give a characterization of $\Aut(\shift{\theta},\Z^d)$ for an arbitrary nontrivial bijective substitution $\theta$. It is easy to see that the shifts and the mapping $\delta: x\mapsto\overline{x}$ are automorphisms of $\shift{\theta}$, so the actual problem is determining whether there exist other kinds of automorphisms. Coven's result in one dimension \cite{Cov71} states that this is not the case:
	\begin{teo}[Coven]
		Let $\Alf=\{0,1\}$ be a two-symbol alphabet. If $\theta:\mathcal{A}\to\mathcal{A}^n$ is a nontrivial (primitive) bijective substitution of constant length $n>1$, then $\Aut(\shift{\theta},\Z)\cong\Z\times(\Z/2\Z)$, with every automorphism being of the form $\sigma^n$ or $\delta\circ\sigma^n$, where $\sigma=\sigma_1$ is the elementary shift action.
	\end{teo}
	Similar results hold for larger alphabets (see e.g. the article by Lema\'nczyk and Metzker \cite{LM88}). Our goal in what follows is to show that these results translate readily to the higher-dimensional case.
	
	For our purpose, studying the action of an automorphism on a fixed point of $\theta$ will be a valuable tool. To begin, we need to introduce some terminology:
	
	\begin{defi}
		Given a $d$-tuple $\vec{u}=(u_1,\dots,u_d)\in\{-1,+1\}^d$, the \textbf{canonical quadrant} associated to $\vec{u}$ is the following subset of $\Z^d$:
		\[Q_{\vec{u}}=u_1\N_0\times u_2\N_0\times\cdots\times u_d\N_0. \]
		We shall refer to any translation $Q_{\vec{u},\vec{k}}\dfn\vec{k}+Q_{\vec{u}}$ of a canonical quadrant as a \textbf{quadrant} and its unique extremal element $\vec{k}$ (with the maximum or minimum possible value on each coordinate) as its \textbf{vertex}.
	\end{defi}
	
	Notice that the $2^d$ quadrants $Q_{\vec{\ind}}=Q_{\vec{\ind},\vec{0}},Q_{(-1,1,\dots),-\vec{e}_1}, \dotsc, Q_{-\vec{\ind},-\vec{\ind}}$ are pairwise disjoint and their union is $\Z^d$; also, $\vec{u}\in Q_{\vec{u}}$. By definition, if two points $x,y\in\mathcal{A}^{\Z^d}$ coincide on a canonical quadrant, then their images under $\theta$ coincide on the same quadrant as well. For a general quadrant, this holds too, although the vertex may change if the quadrant is not canonical.
	
	We can quickly verify that the following holds:
	
	\begin{lem}
		Let $\theta$ be a nontrivial (primitive) bijective substitution. If $x,y\in\shift{\theta}$ coincide on a quadrant, then $\varphi(x)=\varphi(y)$, where $\varphi:\shift{\theta}\epi\Z_{\vec{s}}$ is the encoding factor map from the previous section.
	\end{lem}
	
	\begin{dem}
		This is direct from the uniqueness of the factorization of a point of a one-dimensional substitutive subshift as a concatenation of words of the form $\theta^k(0),\theta^k(1)$, applied to each of the $d$ principal directions.
		
		However, we may also argue as follows: assuming, w.l.o.g. that $x$ and $y$ match on the canonical quadrant $Q_{\vec{\ind}}=\N_0^d$, it is easy to see that $\lim_{k\to\infty}\sigma_{\vec{\ind}}^{h(k)}(x)=\lim_{k\to\infty}\sigma_{\vec{\ind}}^{h(k)}(y)$ for any increasing subsequence $h:\N\to\N_0$ such that $\sigma_{\vec{\ind}}^{h(k)}(x)$ converges, as $\sigma_{\vec{\ind}}^{h(k)}(y)$ coincides with the former in $[-h(k),h(k)]^d$ and thus they are at distance at most $2^{-h(k)}$.
		
		If $\varphi(x)-\varphi(y)=\vec{m}\ne\vec{0}$, it is easy to see that the value $\vec{m}$ remains constant for the aforementioned subsequence:
		\[\varphi(\sigma_{\vec{\ind}}^{h(k)}(x)) - \varphi(\sigma_{\vec{\ind}}^{h(k)}(y)) = (\varphi(x)+h(k)\vec{\ind}) - (\varphi(y)+h(k)\vec{\ind}) = \varphi(x)-\varphi(y)=\vec{m} \]
		and thus, since $\varphi$ is continuous, $\lim_{k\to\infty}\varphi(\sigma_{\vec{\ind}}^{h(k)}(x)) - \varphi(\sigma_{\vec{\ind}}^{h(k)}(y))=\vec{m}$. But also, \[\lim_{k\to\infty}\varphi(\sigma_{\vec{\ind}}^{h(k)}(x)) - \varphi(\sigma_{\vec{\ind}}^{h(k)}(y)) = \lim_{k\to\infty}\varphi(\sigma_{\vec{\ind}}^{h(k)}(x)) - \underbrace{\lim_{k\to\infty}\varphi(\sigma_{\vec{\ind}}^{h(k)}(y))}_{{}=\lim_{k\to\infty}\varphi(\sigma_{\vec{\ind}}^{h(k)}(x))}=\vec{0},\]
		hence $\vec{m}=\vec{0}$, a contradiction.\qed
	\end{dem}
	
	In what follows, we will show that Coven's result for one-dimensional bijective substitutions applies to higher-dimensional substitutive subshifts, i.e. the bijectiveness condition immediately forces the only nontrivial automorphism (up to composition by a shift) to be the relabeling map $\delta$ that swaps the two symbols of the alphabet. The main bulk of the proof of this result lies in the following lemma:
	
	\begin{lem}
		\label{lem:large_scale_relabeling}
		Let $\theta:\mathcal{A}\to\mathcal{A}^S$ be a bijective substitution with nontrivial support $S=[\vec{0},\vec{s}-\vec{\ind}]$ over the alphabet $\mathcal{A}=\{0,1\}$, and suppose $f\in\Aut(\shift{\theta},\Z^d)$ is an automorphism. Then, for any $x\in\shift{\theta}$ there exist $\vec{k},\vec{\ell}\in\Z^d$ and a sufficiently large $m\ge 1$ such that both $x$ and $f(x)$ are concatenations of patterns of the form $\theta^m(0)$ or $\theta^m(1)$ arranged over a translation of a ``grid'' $\vec{s}^m\cdot\Z^d$, and such that the pattern with support $\vec{k}+\vec{p}+S^{(m)}$ (with $\vec{p}\in\vec{s}^m\cdot \Z^d$) in the grid corresponding to $x$ determines uniquely the pattern with support $\vec{\ell}+\vec{p}+S^{(m)}$ in the grid corresponding to $f(x)$.
	\end{lem}
	
	\begin{dem}
		As above, it is a direct consequence of Lemma \ref{lem:codified_system_subst} that for a fixed $m\ge 1$ any point $x\in\shift{\theta}$ is a concatenation of patterns of the form $\theta^m(a),a\in\mathcal{A}$ over a grid given by a translation of $\vec{s}^m\cdot\Z^d$. So we actually are proving the correspondence between these patterns in $x$ and $f(x)$.
		
		%We start by restating notation. Remember that, as stated above, $S^{(m)}$ is the support of the substitution $\theta^m$, while, given any subset $R\subseteq\Z^d$, $R^{\circ m}$ may be thought of as the set of all points from $R$ at distance at least $m$ from the ``border'' of $R$, a discrete analogue to the concept of interior.
		
		By its nature as a sliding block code, any automorphism $f\in\Aut(\shift{\theta},\Z^d)$ has a radius $r\in\N_0$, namely, for any $\vec{k}\in\Z^d$ the symbol $f(x)_{\vec{k}}$ is uniquely determined by the finite pattern $x|_{\vec{k}+[-r\vec{\ind},r\vec{\ind}]}$; thus, for any subset $R\subseteq\Z^d$, $x|_R$ determines uniquely the configuration $f(x)|_{R^{\circ r}}$. From now on, $f$ will be any fixed automorphism and $r$ will be the corresponding radius. Consider then the support $S^{(m)}$ of $\theta^m$; as $S$ was deemed nontrivial, $S^{(m)}$ must be a $d$-dimensional rectangle of edge length at least $2^m$ in any direction, and thus for sufficiently large $m$ (say, $m>\log_2(2r+1)$) the set $(S^{(m)})^{\circ r}$ is nonempty and a $d$-dimensional rectangle of edge length at least $2^m-2r$ in all directions.
		
		By Lemma \ref{lem:codified_system_subst}, there are vectors $\vec{k},\vec{\ell}\in\Z^d$ such that, for any $\vec{p}\in\vec{s}^m\cdot\Z^d$, $x|_{\vec{k}+\vec{p}+S^{(m)}}$ and $f(x)|_{\vec{\ell}+\vec{p}+S^{(m)}}$ are either $\theta^m(0)$ or $\theta^m(1)$; we shall refer to these rectangles as $K_{\vec{p}}\dfn \vec{k}+\vec{p}+S^{(m)}$ and $L_{\vec{p}}\dfn \vec{\ell}+\vec{p}+S^{(m)}$, respectively, for any $\vec{p}\in\vec{s}^m\cdot\Z^d$. Note that, since $S^{(m)}=[\vec{0},\vec{s}^m-\vec{\ind}]$ is a set of representatives for $\Z^d/(\vec{s}^m\cdot\Z^d)$, the rectangles $K_{\vec{p}}$, indexed by all $\vec{p}\in\vec{s}^m\cdot\Z^d$, cover $\Z^d$ completely (and thus the $L_{\vec{p}}$ rectangles do so as well). Since we may replace $\vec{k}$ by any $\vec{k}+\vec{s}^m\cdot\vec{k}'$ (as then the new $K'_{\vec{p}}$ is just the old $K_{\vec{p}+\vec{k}'}$), we may choose $\vec{k}$ in a suitable way such that, for any $\vec{p}\in\vec{s}^m\cdot\Z^d$, $K_{\vec{p}}^{\circ r}$ has nonempty intersection with $L_{\vec{p}}$, say $I_{\vec{p}}\dfn K_{\vec{p}}^{\circ r}\cap L_{\vec{p}}$ (this is because the union of all $L_{\vec{p}}$ is the whole of $\Z^d$; we only need to note that, for a suitable choice of $\vec{k}$, the intersection $I_{\vec{0}}=K_{\vec{0}}^{\circ r}\cap L_{\vec{0}}$ is nonempty, and then use the fact that $K_{\vec{p}}$ and $L_{\vec{p}}$ are translations of $K_{\vec{0}}$ and $L_{\vec{0}}$ by the same vector). It is important to remark that, even though in most arguments we choose $\vec{k}$ and $\vec{\ell}$ from the set $S^{(m)}=[\vec{0},\vec{s}^m-\vec{\ind}]$, as the obvious representatives of the cosets of $\vec{s}^m\cdot\Z^d$, it is not actually necessary to do so, and in particular in this proof $\vec{k}$ and $\vec{\ell}$ may be any two elements from $\Z^d$.
		
		As stated above, since $\theta$ (and thus $\theta^m$) is a bijective substitution, then for any $a,b\in\mathcal{A}$ and any $\vec{q}\in S^{(m)}$ the condition $\theta^m(a)_{\vec{q}}=\theta^m(b)_{\vec{q}}$ implies $a=b$ and thus $\theta^m(a)=\theta^m(b)$. Because of this, whether the pattern $f(x)|_{L_{\vec{p}}}$ is either $\theta^m(0)$ or $\theta^m(1)$ is entirely determined by the subpattern $f(x)|_{I_{\vec{p}}}$ (as $I_{\vec{p}}$ is nonempty), which in turn, as a subpattern of $f(x)|_{K_{\vec{p}}^{\circ r}}$, is entirely determined by $x|_{K_{\vec{p}}}$, which is either $\theta^m(0)$ or $\theta^m(1)$ as well. Thus, for any $\vec{p}\in\vec{s}^m\cdot\Z^d$, $f(x)|_{L_{\vec{p}}}$ depends uniquely on $x|_{K_{\vec{p}}}$, as desired.\qed
	\end{dem}
	
	\begin{cor}
		\label{cor:explicit_morphism}
		Given any $f\in\Aut(\shift{\theta},\Z^d)$ and any $x\in\shift{\theta}$, $f(x)$ is either $\sigma_{\vec{\ell}-\vec{k}}(x)$ or $\delta\circ\sigma_{\vec{\ell}-\vec{k}}(x)$, where $\vec{k}$ and $\vec{\ell}$ are the vectors from the previous lemma.
	\end{cor}
	
	\begin{dem}
		By the previous lemma, $f(x)|_{L_{\vec{p}}}$ is entirely determined by $x|_{K_{\vec{p}}}$ and thus there is a mapping $t:\{0,1\}\to\{0,1\}$ (depending only on the chosen $x$ and the automorphism $f$) such that if $x|_{K_{\vec{p}}}$ is $\theta^m(a)$, then $f(x)|_{L_{\vec{p}}}$ is $\theta^m(t(a))$; the same $t$ applies to all pairs of patterns $x|_{K_{\vec{p}}}$ and $f(x)|_{L_{\vec{p}}}$ for all $\vec{p}\in\vec{s}^m\cdot\Z^d$ due to the Curtis-Hedlund-Lyndon theorem.  If $t(0)=t(1)$, then $f$ sends both $x$ and $\delta(x)$ to the same point, contradicting the bijectiveness of $f$ (since $\theta$ is a primitive substitution and thus has more than one point). Thus, $t(0)\ne t(1)$ and then either $t(a)=a$ or $t(a)=1-a$.
		
		In the first case, if for some $\vec{p}\in\vec{s}^m\cdot\Z^d$ we have $x|_{K_{\vec{p}}}=\theta^m(a)$, then $f(x)|_{L_{\vec{p}}}=\theta^m(t(a))=\theta^m(a)$. This applies to all $\vec{p}\in\vec{s}^m\cdot\Z^d$; since $L_{\vec{p}}=K_{\vec{p}}+(\vec{\ell}-\vec{k})$, from this we see that $f(x)=\sigma_{\vec{\ell}-\vec{k}}(x)$. In the second case, from $x|_{K_{\vec{p}}}=\theta^m(a)$ we deduce that $f(x)|_{L_{\vec{p}}}=\theta^m(t(a))=\overline{\theta^m(a)}=\delta(x)|_{K_{\vec{p}}}$ (since $\theta^{m}(t(x))=\overline{\theta^m(x)}$). Again, this applies to all $\vec{p}\in\vec{s}^m\cdot\Z^d$; thus, $f(x)=\sigma_{\vec{\ell}-\vec{k}}(\delta(x)) = \delta\circ\sigma_{\vec{\ell}-\vec{k}}(x)$.\qed
	\end{dem}
	This almost completes the proof of our first main result:
	\begin{teo}
		For a nontrivial, primitive, bijective substitution $\theta$ on the alphabet $\Alf=\{0,1\}$, $\Aut(\shift{\theta},\Z^d)$ is generated by the shifts and the relabeling map (flip map) $\delta(x)\dfn\overline{x}$, and thus is isomorphic to $\Z^d\times(\Z/2\Z)$.
	\end{teo}
	\begin{dem}
		For any $f\in\Aut(\shift{\theta},\Z^d)$ and any $x\in\shift{\theta}$, $f(x)$ is either $\sigma_{\vec{\ell}-\vec{k}}(x)$ or $\delta\circ \sigma_{\vec{\ell}-\vec{k}}(x)$. Since $f$ commutes with the shift action, $f|_{\Orb_\sigma(x)}$ is either $\sigma_{\vec{\ell}-\vec{k}}|_{\Orb_\sigma(x)}$ or $(\delta\circ \sigma_{\vec{\ell}-\vec{k}})|_{\Orb_\sigma(x)}$.
		
		From the definition of $\shift{\theta}$, there is a finite subset $\Sigma\subset\shift{\theta}$, comprised of periodic points of $\theta_\infty$, such that the union $\bigcup_{x\in\Sigma}\Orb_\sigma(x)$ is dense in $\shift{\theta}$. Thus, for each $x\in\Sigma$, we have that $f|_{\Orb_\sigma(x)}$ is either $\sigma_{\vec{\ell}-\vec{k}}|_{\Orb_\sigma(x)}$ or $(\delta\circ\sigma_{\vec{\ell}-\vec{k}})|_{\Orb_\sigma(x)}$. It is also easy to see that, for all $x\in\Sigma$ the inclusion $\shift{\theta}^\circ\subseteq\overline{\Orb_\sigma(x)}$ holds and thus, since $\sigma_{\vec{\ell}-\vec{k}}|_{\Orb_\sigma(x)}$ and $\delta\circ\sigma_{\vec{\ell}-\vec{k}}|_{\Orb_\sigma(x)}$ differ in the minimal substitutive subshift $\shift{\theta}^\circ$, $f$ cannot be equal to $\sigma_{\vec{\ell}-\vec{k}}$ in an orbit $\Orb_\sigma(x)$ and equal to $\delta\circ\sigma_{\vec{\ell}-\vec{k}}$ in a different orbit $\Orb_\sigma(y)$ for $x\ne y\in\Sigma$. Hence, $f$ either equals $\sigma_{\vec{\ell}-\vec{k}}|_{\Orb_\sigma(x)}$ in all orbits of points of $\Sigma$ or equals $\delta\circ\sigma_{\vec{\ell}-\vec{k}}|_{\Orb_\sigma(x)}$ in all orbits of $\Sigma$. In either case, by density of $\bigcup_{x\in\Sigma}\Orb_\sigma(x)$, $f$ must equal one of these two automorphisms in the whole of $\shift{\theta}$, proving the desired result.\qed
	\end{dem}
	
	Note that in the proof we only used that a point from a substitutive subshift in a two-symbol alphabet is a concatenation of patterns $\theta^m(0)$ and $\theta^m(1)$ and thus the same proof applies for $\shift{\theta}^\circ$ by replacing the set $\Sigma$ with $\Sigma^\circ=\Sigma\cap\shift{\theta}^\circ$, which is nonempty and contains all fixed points of $\theta_\infty$ whose seeds are subpatterns of $\theta^m(a),a\in\mathcal{A}$ for some $m$. Consequently, we state this as a corollary:
	
	\begin{cor}
		For a nontrivial, primitive, bijective substitution $\theta$ on $\Alf=\{0,1\}$, $\Aut(\shift{\theta}^\circ,\Z^d)$ is generated by the shifts and the relabeling map $\delta$, and thus is isomorphic to $\Z^d\times(\Z/2\Z)$.
	\end{cor}
	
	To conclude this section, we shall make some brief remarks regarding the restriction to a two-symbol alphabet, $\mathcal{A}=\{0,1\}$. As noted above, in this special case the definition of bijectivity of a substitution reduces to the condition $\theta(1)=\overline{\theta(0)}$, which results in an explicit description of the only nontrivial (modulo the shifts) automorphism of $\shift{\theta}$, which is the mapping $\delta(x)\dfn\overline{x}$.
	
	In the case of a larger alphabet, the structure of the nontrivial automorphisms might be different, and an automorphism with the same behavior as $\delta$ may not even exist. For instance, if we consider the one-dimensional substitution on the three-letter alphabet $\mathcal{A}=\{1,2,3\}$ given by:
	\[\theta: 1\mapsto123,\qquad 2\mapsto 231,\qquad 3\mapsto312, \]
	then the relabeling map defined as $\varphi_{(1\,2\,3)}(x)_i=\tau(x_i)$, where $\tau=(1\,2\,3)$ is a cyclic permutation of the symbols of $\mathcal{A}$, is a nontrivial element of $\Aut(\shift{\theta},\Z)$ of order $3$. However, it is not clear whether elements of order $2$ do exist in this group, making this mapping the only obvious analogue to $\delta$ satisfying the property $\varphi_{(1\,2\,3)}\circ\theta_\infty=\theta_\infty\circ\varphi_{(1\,2\,3)}$. By changing the form of the substitution, this relabeling map may not even exist, e.g. for:
	\[\vartheta:1\mapsto 123,\qquad 2\mapsto 212,\qquad 3\mapsto331, \]
	the word $331331$ is a subpattern of every point in $\shift{\vartheta}$, and thus $13$ and $333$ are subpatterns of every point as well. But neither $111$ nor $222$ can appear as subwords of a point of $\shift{\vartheta}$, and thus a relabeling map must map $3$ to $3$. But then it has to map the symbols $1$ and $2$ to themselves, to preserve all instances of the word $331331$ and the bijectivity. Thus, nontrivial relabeling maps do not exist.
	
	Obviously, this does not preclude the existence of nontrivial automorphisms of $\shift{\vartheta}$ that are not relabeling maps, i.e. have radius $r$ strictly greater than $0$. However, a cursory look at the proof above shows that most of it does not make explicit usage of the alphabet size. In particular, the proof of Lemma \ref{lem:large_scale_relabeling} carries over without any significant change. In this more general context, this lemma shows that there must be a bijection $\tau:\mathcal{A}\to\mathcal{A}$ such that (using the notation from the proof of Lemma \ref{lem:large_scale_relabeling}) the following holds:
	\[(\exists m\in\N)(\forall\vec{p}\in\vec{s}^m\cdot\Z^d):x|_{K_{\vec{p}}}=\theta^m(a)\iff f(x)|_{L_{\vec{p}}}=\theta^m(\tau(a)),  \]
	and since each $L_{\vec{p}}$ is a translation of the corresponding $K_{\vec{p}}$, we see that, if $\tau_\infty$ is the relabeling map $\mathcal{A}^\Z\to\mathcal{A}^\Z$ induced by $\tau$, then the previous statement can be restated as follows:
	\[(\exists\vec{k}\in\Z^d):f(\theta^m_\infty(x))=\sigma_{\vec{k}}(\theta^m_\infty(\tau_\infty(x))). \]
	So, $f$ behaves similarly to a relabeling map; in particular, this is enough to show that $f^k$ is a shift for sufficiently large $k$, i.e. $f$ has finite order modulo $\Z^d$. We may refine this result even further, by showing that $f$ is indeed the composition of a relabeling map and a shift:
	
	\begin{cor}
		Let $\theta$ be a nontrivial, primitive, bijective substitution on an alphabet $\Alf$ (which can have more than two symbols). For any $f\in\Aut(\shift{\theta},\Z^d)$, there exists a bijection $\tau:\mathcal{A}\to\mathcal{A}$ and a value $\vec{k}\in\Z^d$ such that $f=\sigma_{\vec{k}}\circ\tau_\infty$. Thus, $\Aut(\shift{\theta},\Z^d)$ is isomorphic to a subgroup of $\Z^d\times S_{|\mathcal{A}|}$, where $S_n$ is the symmetric group in $n$ elements.
	\end{cor}
	
	\begin{dem}
		First of all, note that the $m$ from the proof of the result from Lemma \ref{lem:large_scale_relabeling} can be replaced with any $m'>m$ without substantial changes in the proof. This means, in particular, that the following holds (using that $\theta^{m+1}=\theta^m\circ\theta$):
		\begin{align*}
			(\exists\vec{k},\vec{k}'\in\Z^d)(\exists\tau,\tau':\mathcal{A}\to\mathcal{A}):f(\theta^{m+1}_\infty(x)) &=\sigma_{\vec{k}}(\theta^{m}_\infty(\tau_\infty(\theta_\infty(x))))\\
			&=\sigma_{\vec{k}'}(\theta^{m+1}_\infty(\tau'_\infty(x))),
		\end{align*}
 
		and since $\vec{k}\equiv\vec{k}'\pmod{\vec{s}^{m+1}}$, this implies that each pattern $\theta^{m+1}(\tau'(a))$ with $a\in\mathcal{A}$ is a concatenation of the patterns $\theta^m(\tau(b))$, where the $b$ are the corresponding symbols of the pattern $\theta(a)$. But by definition $\theta^{m+1}(\tau'(a))=\theta^m(\theta(\tau'(a)))$, and the mapping $\theta$ is injective; thus, $\theta(\tau'(a))=\tau(\theta(a))$, i.e. the relabeling $\tau$ must send patterns of the form $\theta(b)$ to other patterns of the form $\theta(b')$.
		
		By replacing $\theta$ with a suitable power, we may assume that for the bottom left corner $\vec{0}$ of the support $S$ the equality $\theta(a)_{\vec{0}}=a$ holds. Thus, $\theta(\tau'(a))$ has $\tau'(a)$ in this position, while $\tau(\theta(a))$ has $\tau(a)$ in the same position, i.e. $\tau(a)=\tau'(a)$. As this applies to any symbol $a$, we conclude that $\tau=\tau'$ and that $\tau$ and $\theta$ commute, i.e. $\theta_\infty\circ\tau_\infty=\tau_\infty\circ\theta_\infty$ as mappings $\mathcal{A}^{\Z^d}\to\mathcal{A}^{\Z^d}$. Applying this result to the identity with $f$ above, we conclude that:
		\[(\exists\vec{k}\in\Z^d):f(\theta_\infty^m(x))=\sigma_{\vec{k}}\circ\tau_\infty(\theta_\infty^m(x)), \]
		and since every point $x\in\shift{\theta}$ is a shift of a point of the form $\theta_\infty^m(x)$ by Lemma \ref{lem:codified_system_subst}, this and the Curtis-Hedlund-Lyndon theorem show that $f=\sigma_{\vec{k}}\circ\tau_\infty$, the desired result. Since $\tau_\infty$ is entirely determined by a bijection from the finite set $\mathcal{A}$ to itself, and it is obvious that a relabeling map is shift-commuting by definition, we can identify $\Aut(\shift{\theta},\Z^d)$ with a subgroup of $\Z^d\times S_{|\mathcal{A}|}$.\qed
	\end{dem}
	
	Note that this proof also provides a necessary condition for a bijection $\tau:\mathcal{A}\to\mathcal{A}$ to induce a relabeling map $\tau_\infty\in\Aut(\shift{\theta},\Z^d)$; namely, that $\theta_\infty\circ\tau_\infty=\tau_\infty\circ\theta_\infty$. By compactness, this condition is also sufficient, providing an explicit description of the group $\Aut(\shift{\theta},\Z^d)$ in terms of the patterns $\theta(a),a\in\mathcal{A}$.
	
	\section{Extended symmetries and bijective substitutions}
	
	Our next goal is to obtain generalizations of the previous result in the domain of extended symmetries. These are a generalization of automorphisms, which introduce an additional degree of flexibility by allowing, besides the standard local transformations given by a sliding block code, to ``deform'' the underlying $\Z^d$ lattice, by rotation, reflection, shear or other effects of a geometric nature. This additional degree of freedom is captured by a group automorphism of $\Z^d$, i.e. an element\footnote{Remember that $\GL_d(\Z)$ is the set of all invertible matrices with integer entries whose inverses are also matrices with integer entries (namely, all matrices with integer entries that satisfy the condition $\det(A)=\pm 1$). Any matrix of this kind induces a bijective linear transformation $T_A:\Z^d\to\Z^d,\vec{p}\mapsto A\vec{p}$ and vice versa. We may generalize the previous definition to any group $G$ by replacing $\GL_d(\Z)$ with $\Aut(G)$, the set of group automorphisms of $G$. However, as we shall be primarily concerned with $G=\Z^d$, the restricted definition is enough for our purposes.} of $\GL_d(\Z)$.
	
	The basic premises of the theory of extended symmetries of subshifts may be studied in \cite{BRY2018}.
	
	\begin{defi}
		Let $X\subseteq\mathcal{A}_X^{\Z^d},Y\subseteq\mathcal{A}_Y^{\Z^d}$ be two $\Z^d$-subshifts. Given a $\Z$-invertible matrix with integer entries $A\in\GL_d(\Z)$, we call a continuous mapping $f:X\to Y$ an $A$-\textbf{morphism} if the following equality holds:
		\[(\forall \vec{p}\in\Z^d):f\circ\sigma_{\vec{p}}=\sigma_{A\vec{p}}\circ f.  \]
		An \textbf{extended symmetry} is a bijective $A$-morphism from $X$ to itself, associated to some $A\in\GL_d(\Z)$. We shall denote the set of all extended symmetries as $\Sym(X,\Z^d)$. This is a group under composition.
	\end{defi}
	Under our standard hypothesis (namely, a faithful shift action) the matrix $A_f$ associated to an extended symmetry $f$ is uniquely determined and thus there is an obvious mapping $\psi:\Sym(X,\Z^d)\to\GL_d(\Z), f\mapsto A_f$. It is also easy to see that $\psi$ is a group morphism:
	\[f\circ g\circ \sigma_{\vec{p}} = f\circ\sigma_{\psi(g)\vec{p}}\circ g = \sigma_{\psi(f)(\psi(g)\vec{p})}\circ f\circ g,   \]
	consequently, $\psi(f\circ g)=\psi(f)\psi(g)$. Evidently, $\psi(f)=I_d$ (the identity matrix) if and only if $f$ is a traditional automorphism of $X$, i.e. $\ker(\psi)=\Aut(X,\Z^d)$. This implies that the quotient group $\Sym(X,\Z^d)/\Aut(X,\Z^d)$ is isomorphic to a subgroup of $\GL_d(\Z)$, and thus, determining the nature of the latter group as a subgroup of $\GL_d(\Z)$ is a very useful tool to describe $\Sym(X,\Z^d)$.
	
	Due to their ``almost shift-commuting'' nature, there is a version of the Curtis-Hedlund-Lyndon theorem for extended symmetries (and, in general, for $A$-morphisms), which implies that an extended symmetry is a composition of a local map (in the sense of the classical Curtis-Hedlund-Lyndon theorem) and a lattice transformation given by the matrix $A$. The result, proved in \cite{BRY2018} is as follows:
	
	\begin{teo}[Generalized Curtis-Hedlund-Lyndon theorem] \label{teo:curtishedlundlyndon}
		Let $f:X\to X$ be an extended symmetry from $\Sym(X,\Z^d)$. Then, there is a finite subset $U\subset\Z^d$ and a function $F:\mathcal{A}^U\to\mathcal{A}$ such that the following equality holds for all $\vec{s}\in\Z^d$:
		\[f(x)_{A\vec{s}} = F(x|_{\vec{s}+U}), \]
		in which we identify, as usual, a pattern with support $U$ with any of its translations.
	\end{teo}
	This, as is the case for automorphisms, allows us to show that whether two points match on a ``large'' set $R\subseteq\Z^d$, their images under an extended symmetry $f$ match as well on a large set, which depends on $f$ and $R$. More precisely, if we suppose w.l.o.g. that the support $U$ (as defined in the theorem) of the symmetry $f$ is of the form $[-r\vec{\ind},r\vec{\ind}]$, then:
	\[x|_R = y|_R\implies f(x)|_{\psi(f)[R^{\circ r}]}=f(y)|_{\psi(f)[R^{\circ r}]}. \]
	In particular, if $R$ is a half-space, the set $\psi(f)[R^{\circ r}]$ is a half-space as well.
	
	Our goal is to characterize the group $\Sym(X,\Z^d)/\Aut(X,\Z^d)$ when $X=\shift{\theta}$, $\theta$ being a nontrivial bijective substitution, and then characterize the extended symmetry group explicitly whenever possible. We shall see that $\shift{\theta}$, by construction, has some distinguished points with \textbf{fractures}, that is, they are comprised of subconfigurations of points from the minimal substitutive subshift $\shift{\theta}^\circ$ ``glued together'' in an independent way, and that any extended symmetry has to preserve these points with fractures. In particular, by analyzing these points adequately we can deduce strong restrictions on the matrices $\psi(f)\in\GL_d(\Z)$ for any $f\in\Sym(X,\Z^d)$, as the ``shape'' of the fractures determined by a certain subset of $\Z^d$ must be preserved by the matrix $\psi(f)$. Our choice of adequate points to show this result is given by the following lemma:

	\begin{lem}\label{lem:points_of_contradiction}
		Given any nontrivial bijective primitive substitution $\theta:\Alf\to\Alf$ over a two-symbol alphabet $\Alf=\{0,1\}$, there exist $x,y\in\shift{\theta}$ such that $x|_{Q_{\vec{\ind}}^c}=y|_{Q_{\vec{\ind}}^c}$ but $x|_{Q_{\vec{\ind}}}=\overline{y}|_{Q_{\vec{\ind}}}$, where $Q_{\vec{\ind}}=\{(n_1,\dots,n_d)\in\Z^d: (\forall 1\le i\le d):n_i\ge 0 \}$ is the canonical quadrant containing $\vec{\ind}$. 
	\end{lem}
	
	\begin{dem}
		Without loss of generality, we may assume that the pattern $\theta^m(a)$ has the symbol $a$ on all $2^d$ corners for all $m$ (by replacing $\theta$ by a power $\theta^k$ if needed); this implies that there are fixed configurations for $\theta$ over $\N^d$ with seeds $0$ and $1$. Thus, taking any fixed point of $\theta$ over $\Z^d$ with seed $P$ and changing the values of $P_{\vec{0}}$ we obtain two valid points $x,y\in\shift{\theta}$ that differ only on the positive quadrant $Q_{\vec{\ind}}$.
		
		By the nature of a bijective substitution, since $x$ and $y$ are fixed points of $\theta_\infty$, $x_{\vec{s}}=y_{\vec{s}}$ for any $\vec{s}\in Q_{\vec{\ind}}$ would imply $x|_{Q_{\vec{\ind}}}=y|_{Q_{\vec{\ind}}}$, which we know is not the case; thus, $x|_{Q_{\vec{\ind}}}=\overline{y}|_{Q_{\vec{\ind}}}$.\qed
	\end{dem}
	\begin{nota}
		A similar result holds for general (finite) alphabets: we can find two points $x,y$ such that $x|_{Q_{\vec{\ind}}^c}=y|_{Q_{\vec{\ind}}^c}$ but $x|_{\vec{s}}=y|_{\vec{s}}$ for all positions $\vec{s}\in Q_{\vec{\ind}}$.
	\end{nota}
	
	\begin{figure}[ht]
		\centering
		\includegraphics[scale=0.7]{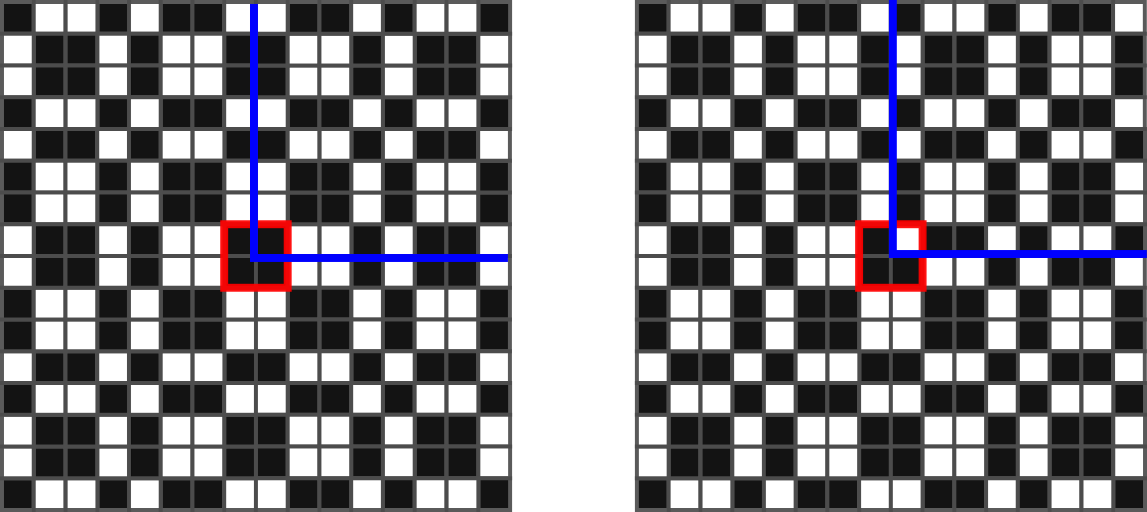}
		\caption{Two points from the two-dimensional Thue-Morse substitutive subshift $\shift{\theta_{TM}}$ matching in all but one quadrant. The seed of each point corresponds to the central $2\times 2$ pattern.}
		\label{fig:contradiction_points}
	\end{figure}
	
	A pair of points satisfying the previous lemma is displayed in Figure \ref{fig:contradiction_points}, for the Thue-Morse substitution. It is easy to see that a similar argument works for any quadrant, changing the form of the substitution if so required: specifically, we may replace $\theta$ by an iterate $\theta^m$ such that $\theta^m(a)$ has the symbol $a$ on all corners, in particular, the corner $\vec{v}$ of $S^{(m)}$ that corresponds to the vertex of this quadrant (e.g. for $Q_{-\vec{\ind}}$, the quadrant containing $(-1,\dots,-1)$, we have $\vec{v}=\vec{s}^d$). The rest of the argument from Lemma \ref{lem:points_of_contradiction} applies without modifications to show that changing the (unique) symbol of the seed located in a specific quadrant changes the symbols of the whole quadrant, without affecting the remaining symbols of the other quadrants.
	
	By the generalized Curtis-Hedlund-Lyndon theorem, we verify that the image of the aforementioned two points by any extended symmetry $f$ matches along a ``large'' set of the form $\psi(f)[(Q_{\vec{\ind}}^{c})^{\circ r}]$. We shall use this to prove that, unless the image of a quadrant $Q$ by $\psi(f)$ is itself a quadrant, the restriction of $x$ to $Q$ determines $f(x)$ not only in $\psi(f)([(Q^c)^{\circ r}]$, but in $\psi(f)([Q^{\circ r}]$ as well, and from this we later infer that $x|_{Q^c}$ determines $f(x)$ in the whole plane; thus, the existence of two distinct points that match in $Q^c$ contradicts the bijectivity of $f\in\Sym(X,\Z^d)$.
	
	For this purpose, we need first to determine what kind of ``shearing'' is allowed for an extended symmetry in a substitutive subshift. We introduce the following result:
	
	\begin{lem} \label{lem:fitting_quadrants}
		Let $A\in\GL_d(\Z)$ be a matrix with integer coefficients, invertible over $\Z$. Then $AQ_{\vec{\ind}}$ cannot contain two distinct (canonical) quadrants.
	\end{lem}
	
	\begin{dem}
		Notice that $Q_{\vec{\ind}}$ is the set of all nonnegative linear combinations (with coefficients in $\Z$) of the vectors of the canonical basis $\vec{e}_1=(1,0,0,\dotsc),\vec{e}_2=(0,1,0,\dotsc),\dotsc$. Thus, its image under $A$ is the set of all nonnegative integer linear combinations of the columns of $A$.
		
		Suppose $AQ_{\vec{\ind}}$ contains two distinct quadrants, $Q'$ and $Q''$. Note that $Q_{\vec{\ind}}$ is the intersection of the real first quadrant $(\R^+_0)^d$ with $\Z^d$; since $A$ is a matrix from $\GL_d(\Z)$, this means that $AQ_{\vec{\ind}}$ is the intersection of $A(\R^+_0)^d$ (which is a convex set) with $\Z^d$. Thus, any point with integer coordinates which is a convex combination of two points from $AQ_{\vec{\ind}}$ belongs to $AQ_{\vec{\ind}}$ as well. By this convexity argument, we may assume those two quadrants $Q',Q''$ to be adjacent, sharing a ``face'', i.e. there is a subset $H=H_{I,j}\subseteq\Z^d$, given by a set of indices $I\subseteq\{1,\dots,d\}$ and an additional index $j\in\{1,\dots,d\}\setminus I$, which is of the form:
		\begin{align*}
			H=\{(m_1,\dots,m_d)\in\Z^d:\,\, & m_j=0\wedge\phantom{-} \\
			& (\forall i\in I):m_i\ge 0\wedge\phantom{-}\\
			& (\forall i'\in \{1,\dots, d\}\setminus(I\cup\{j\})):m_{i'}<0  \},
		\end{align*}
		such that, using the fact that the element $\vec{e}_j$ from the canonical basis is orthogonal to all elements of $H=H_{I,j}$, we have the following decomposition:
		\begin{align*}
			Q'&=\{\vec{m}\in\Z^d:\vec{m}=\vec{m}_0+\lambda\vec{e}_j, \vec{m}_0\in H,\lambda\ge 0 \},\\
			Q''&=\{\vec{m}\in\Z^d:\vec{m}=\vec{m}_0+\lambda\vec{e}_j, \vec{m}_0\in H,\lambda< 0 \}.
		\end{align*}
		
		Thus, taking some fixed $\vec{m}_0\in H$, we see that for any $n\in\Z$ the element $\vec{m}_0+n\vec{e}_j$ belongs to $AQ_{\vec{\ind}}$. Since $A\in\GL_d(\Z)$, the columns of $A$ are a basis for $\Z^d$ and thus there are uniquely determined coefficients $\lambda_1,\dots,\lambda_d,\mu_1,\dots,\mu_d\in\N$ such that:
		\begin{align*}
			\vec{m}_0 &= \sum_{i=1}^{d} \lambda_i(A\vec{e}_i), \\
			\vec{e}_j &= \sum_{i=1}^{d} \mu_i(A\vec{e}_i), \\
			\vec{m}_0 + n\vec{e}_j &= \sum_{i=1}^{d} (\lambda_i+n\mu_i)(A\vec{e}_i).
		\end{align*}
		Note that, since $Q_{\vec{\ind}}$ is comprised of elements of $\Z^d$ with nonnegative coordinates, the linearity of multiplication by $A$ implies that the coefficients $\lambda_1,\dots,\lambda_d$, being uniquely determined, must be nonnegative; that is, if $\vec{m}\in Q_{\vec{\ind}}$, the coefficients of $\vec{m}$ in the canonical basis carry on to the coefficients of $A\vec{m}$ in the basis $\{A\vec{e}_1,\dots,A\vec{e}_d\}$, the set of columns of $A$, and thus keep their corresponding sign. Because of this, as for all values of $n\in\Z$ the vectors $\vec{m}_0+n\vec{e}_j$ belong to $AQ_{\vec{\ind}}$, the corresponding coefficients $\lambda_i+n\mu_i$ must be nonnegative for any value of $n\in\Z$.
		
		Since $\vec{e}_j\ne\vec{0}$, at least one of the $\mu_i$ is nonzero; thus, choosing $n$ adequately we may force $\lambda_i+n\mu_i$ to be negative. This contradicts the observation above about the positivity of coefficients of elements of $AQ_{\vec{\ind}}$ in the base of columns of $A$; so, this shows that $AQ_{\vec{\ind}}$ cannot contain two disjoint quadrants.\qed
	\end{dem}
	
	\begin{nota}
		It is important to stress that, once again, the previous argument does not depend on the chosen quadrant or whether it is canonical or not, and we use $Q_{\vec{\ind}}$ for ease of description. More precisely, we may always multiply $A$ by a change of base matrix of the form $P=[u_1\vec{e}_{\sigma(1)}\mid\cdots\mid u_d\vec{e}_{\sigma(d)}]$ to swap the quadrant $Q_{(u_1,\dots,u_d)}$ (with $u_1,\dots,u_d\in\{-1,+1\}$) with $Q_{\vec{\ind}}$. The same argument as above applies for the matrix $AP$ and thus for the matrix $A$ as multiplication by $P$ only swaps the positions of all the quadrants.
	\end{nota}
	
	\begin{figure}[h!]
		\centering
		\begin{tikzpicture}[scale=0.6]
		\fill[color=red!10!white] (3,3) -- (0,0) -- (3,0) -- cycle;
		\foreach \i in {-5,...,5} {
			\foreach \j in {-5,...,5} {
				\draw[color=gray!50!white] (\i/2,\j/2) circle (0.75mm);
			}	
		}
		\foreach \i in {0,...,5} {
			\foreach \j in {0,...,\i} {
				\fill (\i/2,\j/2) circle (0.75mm);
			}	
		}
		\draw[-latex] (0,-3) -- (0,3);
		\draw[-latex] (-3,0) -- (3,0);
		\draw[color=red,-latex] (0,0) -- (0.5,0);
		\draw[color=red,-latex] (0,0) -- (0.5,0.5);
		\draw[color=red,dashed] (3,3) -- (0,0) -- (3,0);
		\node at (0,-4) {$A=\begin{bmatrix}
			1 & 1 \\ 0 & 1
			\end{bmatrix}$};
		\end{tikzpicture}\quad
		\begin{tikzpicture}[scale=0.6]
		\fill[color=red!10!white] (-3,3) -- (0,0) -- (3,0) -- (3,3) -- cycle;
		\foreach \i in {-5,...,5} {
			\foreach \j in {-5,...,5} {
				\draw[color=gray!50!white] (\i/2,\j/2) circle (0.75mm);
			}	
		}
		\foreach \j in {0,...,5} {
			\foreach \i in {-\j,...,5} {
				\fill (\i/2,\j/2) circle (0.75mm);
			}	
		}
		\draw[-latex] (0,-3) -- (0,3);
		\draw[-latex] (-3,0) -- (3,0);
		\draw[color=red,-latex] (0,0) -- (0.5,0);
		\draw[color=red,-latex] (0,0) -- (-0.5,0.5);
		\draw[color=red,dashed] (-3,3) -- (0,0) -- (3,0);
		\node at (0,-4) {$A=\begin{bmatrix}
			1 & -1 \\ 0 & 1
			\end{bmatrix}$};
		\end{tikzpicture}\quad
		\begin{tikzpicture}[scale=0.6]
		\fill[color=red!10!white] (-1.5,3) -- (0,0) -- (3,-3) -- (3,3) -- cycle;
		\foreach \i in {-5,...,5} {
			\foreach \j in {-5,...,5} {
				\draw[color=gray!50!white] (\i/2,\j/2) circle (0.75mm);
			}	
		}
		\foreach \i in {0,...,5} {
			\foreach \j in {0,...,5} {
				\fill (\i/2,\j/2) circle (0.75mm);
			}
		}
		\foreach \i in {1,...,5} {
			\foreach \j in {\i,...,5} {
				\fill (\j/2,-\i/2) circle (0.75mm);
			}
		}
		\foreach \i in {2,...,5} {
			\fill (-0.5,\i/2) circle (0.75mm);
		}
		\foreach \i in {4,5} {
			\fill (-1,\i/2) circle (0.75mm);
		}
		\draw[-latex] (0,-3) -- (0,3);
		\draw[-latex] (-3,0) -- (3,0);
		\draw[color=red,-latex] (0,0) -- (-0.5,1);
		\draw[color=red,-latex] (0,0) -- (0.5,-0.5);
		\draw[color=red,dashed] (-1.5,3) -- (0,0) -- (3,-3);
		\node at (0,-4) {$A=\begin{bmatrix}
			1 & -1 \\ -1 & 2
			\end{bmatrix}$};
		\end{tikzpicture}
		\caption{Some of the possibilities for the image of $Q_{\vec{\ind}}$ by $A$. The associated matrix $A$ is given below each diagram.}
		\label{fig:quadrant_image}
	\end{figure}
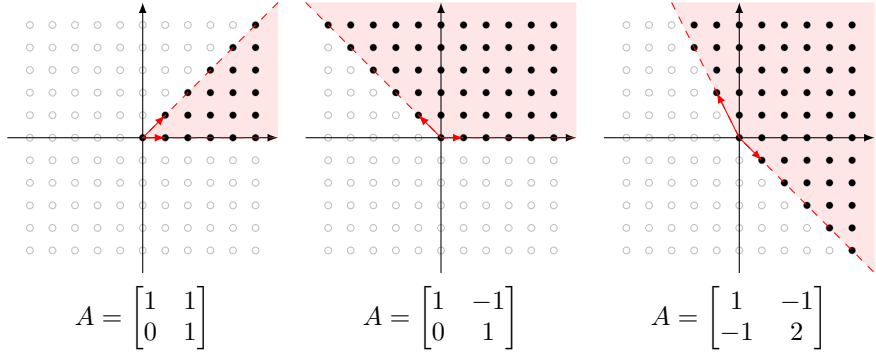
	
	By the previous lemma, there is at most one quadrant $Q_{\vec{v}}$ that is entirely contained in $AQ_{\vec{\ind}}$, meaning that every other quadrant has points from the complement of $AQ_{\vec{\ind}}$. The argument below will use those points and the rigidity of bijective substitutions: any symbol on a fixed position of a pattern $\theta^m(a)$ determines the whole pattern. Our purpose is to exploit this to show that any symbol located at a position from $AQ_{\vec{\ind}}$ is determined by some pattern with support contained in $Q_{\vec{\ind}}^c$, unless the matrix $A$ is of a very specific type, namely, a matrix that does not ``shear'' the underlying $\Z^d$ lattice. The reason behind this lies behind the following result:
	
	\begin{lem}
		\label{lem:finite_lines}
		Let $A\in\GL_d(\Z)$ be a matrix for which $AQ_{\vec{\ind}}$ does not contain a quadrant with vertex $\vec{0}$ as a subset. Then, for any point $\vec{p}\in AQ_{\vec{\ind}}$, one of the affine subspaces $\vec{p}+\Z\vec{e}_i$, for some $i\in\{1,\dots,d\}$, has finite intersection with $AQ_{\vec{\ind}}$.
	\end{lem}
	
	\begin{dem}
		As stated in the previous lemma, the set $AQ_{\vec{\ind}}$ has a sort of ``convexity property'', so that given two points $\vec{q}_1,\vec{q}_2\in AQ_{\vec{\ind}}$, any point $\vec{q}\in\Z^d$ that is a rational convex combination of $\vec{q}_1$ and $\vec{q}_2$ belongs to $AQ_{\vec{\ind}}$ as well.
		
		If each intersection $\psi(f)[Q_{\vec{\ind}}]\cap (\vec{p}+\Z\vec{e}_i)$ were infinite, then by this convexity property, $AQ_{\vec{\ind}}$ must contain, for all indices $j\in\{1,\dots,d\}$, either $\vec{p}+\N_0\vec{e}_i$ or $\vec{p}-\N_0\vec{e}_i$ as a subset. Again by convexity, this implies that $\psi(f)[Q_{\vec{\ind}}]$ contains a quadrant with vertex $\vec{p}$. Since $AQ_{\vec{\ind}}$ is the intersection of $\Z^d$ with a closed, convex subset of $\R^d$ that contains $\vec{0}$ and a quadrant with vertex $\vec{p}$, then $AQ_{\vec{\ind}}$ must contain a quadrant with vertex $\vec{0}$ as a subset, a contradiction with the hypothesis.\qed
	\end{dem}
	
	\begin{nota}
		The same convexity argument as above shows that, if $AQ_{\vec{\ind}}$ (or any set that is the intersection of a closed, convex subset of $\R^d$ and $\Z^d$) contains a quadrant $Q'$ as a subset, then for any $\vec{p}\in AQ_{\vec{\ind}}$ this set contains also a quadrant with vertex $\vec{p}$ as a subset, which is a translation of $Q'$. In particular, the argument above works for any quadrant, not just for $Q_{\vec{\ind}}$.
	\end{nota}
	
	\begin{lem}
		\label{lem:quadrant_determinism}
		Let $\theta$ be a nontrivial, primitive, bijective $d$-dimensional substitution on a finite alphabet $\Alf$ and let $f\in\Sym(\shift{\theta},\Z^d)$ be an extended symmetry. If $\psi(f)Q_{\vec{\ind}}$ is a strict subset of a quadrant $Q$ with vertex $\vec{0}$, then the configuration $x|_{Q_{\vec{\ind}}^c}$ determines $f(x)$ completely.
	\end{lem}
	
	\begin{dem}
		The hypothesis on the matrix $\psi(f)$ shows us by Lemma \ref{lem:finite_lines} that, given any point $\vec{q}\in\psi(f)Q_{\vec{\ind}}$, there exists some direction $\vec{e}_j$ parallel to the coordinate axes for which the intersection $\psi(f)Q_{\vec{\ind}}\cap(\vec{q}+\Z\vec{e}_j)$ is finite. This implies that any rectangle containing $\vec{q}$ and with sufficiently large edge lengths has to contain points outside of $\psi(f)Q_{\vec{\ind}}$, which is the basis of the argument below.
		
		Since $f$ is an extended symmetry, it has an associated radius $r$ so that $f(x)_{A\vec{k}}$ depends only on $[\vec{k}-r\vec{\ind},\vec{k}+r\vec{\ind}]$. Take the least such $r$ and consider the set $(\psi(f)Q_{\vec{\ind}}^c)^{\circ r}$; this set is contained in the complement of $\psi(f)Q_{\vec{\ind}}+[-r\vec{\ind},r\vec{\ind}]$, which is itself contained in some translate of $\psi(f)Q_{\vec{\ind}}$, and thus Lemma \ref{lem:finite_lines} applies to $C^*\dfn \psi(f)Q_{\vec{\ind}}+[-r\vec{\ind},r\vec{\ind}]$.		
		
		Choose a fixed vector $\vec{q}\in C^*$, and let $\vec{e}_j$ be the element from the basis of $\Z^d$ for which the intersection $C^*\cap(\vec{q}+\Z\vec{e}_j)$ is finite. Then, there is a value $h$ such that the intersection $(\vec{q}+\Z\vec{e}_j)\cap C^*$ is contained in the square $[\vec{q}-h\vec{\ind},\vec{q}+h\vec{\ind}]$. Hence, any rectangle with edge length at least $2(h+1)$ in the direction $\vec{e}_j$ must necessarily contain points from the complement of $C^*$, and thus has nonempty intersection with $(\psi(f)Q_{\vec{\ind}}^c)^{\circ r}$. Thus, by choosing $m$ such that the least dimension of $S^{(m)}=[\vec{0},\vec{s}^m-\vec{\ind}]$ is sufficiently bigger than this value $h$ (say, for instance, $m>1+\log_2(h)$) we ensure that the intersection $(\vec{q}+\Z\vec{e}_i)\cap C^*$ has smaller cardinality than the intersection $(\vec{q}+\Z\vec{e}_i)\cap R$, where $R$ is any translate of the rectangle $S^{(m)}$ that contains $\vec{q}$. Thus, $R\cap (\psi(f)Q_{\vec{\ind}}^c)^{\circ r}$ must be nonempty. We visualize this situation in Figure \ref{fig:large_rectangle_overlap}.
		
		\begin{figure}
			\centering
			\begin{tikzpicture}
			\fill[color=red!20!white] (6,6) -- (30/6.5,6) -- (0,0.5) -- (0,0) -- (3,0) -- (6,0.5) -- cycle;
			\fill[color=red!40!white] (6,6) -- (5,6) --  (0,0) -- (6,1) -- cycle;
			\draw[-latex] (0,0) -- (6,0);
			\draw[-latex] (0,0) -- (0,6);
			\foreach \i in {0.4,1.4,...,6} {
				\draw[color=blue!30!white] (\i,0) -- (\i,6);	
			}
			\foreach \i in {0.1,1.6,...,6} {
				\draw[color=blue!30!white] (0,\i) -- (6,\i);
			}
			\foreach \i in {0.4,2.4,...,6} {
				\draw[color=blue!70!white,line width = 1mm] (\i,0) -- (\i,6);	
			}
			\foreach \i in {1.6,4.6,...,6} {
				\draw[color=blue!70!white,line width = 1mm] (0,\i) -- (6,\i);
			}
			\draw[color=red,dashed] (6,1) -- (0,0) -- (5,6);
			\draw[color=red!50!black,dashed] (6,0.5) -- (3,0);
			\draw[color=red!50!black,dashed] (30/6.5,6) -- (0,0.5);
			\fill (2.7,4.3) circle (0.5mm);
			\fill (4.1,1.9) circle (0.5mm);
			\node at (2.7,4.3) [below] {$\vec{p}$};
			\node at (4.1,1.9) [above] {$\vec{q}$};
			\node at (1,5) {$(\psi(f)Q_{\vec{\ind}}^c)^{\circ r}$};
			\node at (5,3) {$C^*$};
			\end{tikzpicture}
			\caption{Example of how the grid decomposition, for sufficiently large $m$, determines rectangles that overlap both the outer region $(\psi(f)Q_{\vec{\ind}}^c)^{\circ r}$ and any desired point from the inner region $C^*=\psi(f)Q_{\vec{\ind}}+[-r\vec{\ind},r\vec{\ind}]$. The symbols from the outer region are entirely determined by the configuration $x|_{Q_{\vec{\ind}^c}}$ due to the generalized Curtis-Hedlund-Lyndon theorem; those from the inner region depend on the points from the outer region because of the bijectivity of $\theta$. }
			\label{fig:large_rectangle_overlap}
		\end{figure}
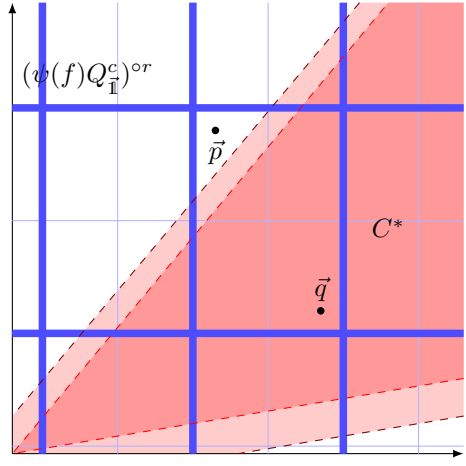
		
		By the encoding from Lemma \ref{lem:codified_system_subst}, we may represent $\Z^2$ as the disjoint union of translates of $S^{(m)}=[\vec{0},\vec{s}^m-\vec{\ind}]$, such that the restriction of $f(x)$ to any of these rectangles is either $\theta^m(0)$ or $\theta^m(1)$. Since the substitution is bijective, knowing a single symbol at any position in the pattern in one of these rectangles is enough to determine whether this pattern is $\theta^m(0)$ or $\theta^m(1)$. Thus, $f(x)_{\vec{q}}$ is entirely determined by any other $f(x)_{\vec{p}}$ where $\vec{p}$ shares the same rectangle $R$. In particular, we may choose $\vec{p}$ from the nonempty intersection $R\cap{(\psi(f)Q_{\vec{\ind}}^c)^{\circ r}}$. But the symbol at coordinate $\vec{p}$ is entirely determined by $x|_{Q_{\vec{\ind}}^c}$ due to the generalized Curtis-Hedlund-Lyndon theorem. Thus, $f(x)_{\vec{q}}$ is already determined by $x|_{Q_{\vec{\ind}}^c}$. This argument applies to any $\vec{q}\in C^*$, proving the desired result.\qed
	\end{dem}
	
	\begin{nota}
		The same proof as above applies to the complement of any quadrant, without regard to its vertex. Thus, if for some quadrant $Q$ the set $\psi(f)[Q]$ is a strict subset of a quadrant with the same vertex, then the same argument applies, showing that we may entirely disregard a quadrant in the preimage $x$ to determine its image $f(x)$.
	\end{nota}
	
	\begin{cor}
		\label{cor:quadrant_permutations}
		Let $\theta$ be a nontrivial, primitive, bijective substitution over a finite alphabet $\Alf$ and let $f\in\Sym(\shift{\theta},\Z^d)$ be an extended symmetry. Then, $\psi(f)$ must send quadrants to quadrants and thus it must induce a permutation on the set of one-dimensional subspaces $\{\Z\vec{e}_1,\Z\vec{e}_2,\dotsc,\Z\vec{e}_d \}$.
	\end{cor}
	
	\begin{dem}
		Suppose $\psi(f)Q$ is not a quadrant for some quadrant $Q$. Then either $\psi(f)Q$ strictly contains a quadrant $Q'$ (w.l.o.g. with vertex $\vec{0}$) or not; the second case falls under the hypothesis of the previous lemma and thus $f$ sends two points that match in the three complementary quadrants of $Q^c$ to the same point, if they exist. Those points indeed do exist, as shown by Lemma \ref{lem:points_of_contradiction}, raising a contradiction with the injectiveness of $f$.
		
		If $\psi(f)[Q]$ strictly contains a quadrant $Q'$ (with vertex $\vec{0}$) then $\psi(f)^{-1}Q'$ is strictly contained in $Q$ and contains $\vec{0}$; thus, the set $\psi(f)^{-1}Q'$ does not contain a quadrant $Q''$ as a subset. But then, since $f^{-1}$ is an extended symmetry itself and $\psi(f^{-1})=\psi(f)^{-1}$, the mapping $f^{-1}$ satisfies the hypothesis of Lemma \ref{lem:quadrant_determinism}, and thus we get to the same contradiction regarding the injectiveness of $f^{-1}$.
		
		\begin{figure}[ht]
			\centering
			\begin{tikzpicture}[scale=0.7]
			\fill[color=red!10!white] (-1.5,3) -- (0,0) -- (3,-3) -- (3,3) -- cycle;
			\foreach \i in {-5,...,5} {
				\foreach \j in {-5,...,5} {
					\draw[color=gray!50!white] (\i/2,\j/2) circle (0.75mm);
				}	
			}
			\foreach \i in {0,...,5} {
				\foreach \j in {0,...,5} {
					\fill (\i/2,\j/2) circle (0.75mm);
				}
			}
			\foreach \i in {1,...,5} {
				\foreach \j in {\i,...,5} {
					\fill (\j/2,-\i/2) circle (0.75mm);
				}
			}
			\foreach \i in {2,...,5} {
				\fill (-0.5,\i/2) circle (0.75mm);
			}
			\foreach \i in {4,5} {
				\fill (-1,\i/2) circle (0.75mm);
			}
			\draw[-latex] (0,-3) -- (0,3);
			\draw[-latex] (-3,0) -- (3,0);
			\draw[color=red,-latex] (0,0) -- (-0.5,1);
			\draw[color=red,-latex] (0,0) -- (0.5,-0.5);
			\draw[color=red,dashed] (-1.5,3) -- (0,0) -- (3,-3);
			\draw[-latex] (4,0) -- (4.1,0.1) -- (4.3,-0.1) -- (4.4,0) -- (5.2,0);
			\node at (6,0) {\vphantom{.}};
			\node at (-1.5,2.5) [left] {$A[Q_{\vec{\ind}}]$};
			\end{tikzpicture}
			\begin{tikzpicture}[scale=0.7]
			\fill[color=red!10!white] (1.5,3) -- (0,0) -- (3,3) -- cycle;
			\foreach \i in {-5,...,5} {
				\foreach \j in {-5,...,5} {
					\draw[color=gray!50!white] (\i/2,\j/2) circle (0.75mm);
				}	
			}
			\foreach \i in {0,...,2} {
				\foreach \j in {0,...,\i} {
					\fill (\i/2,\i/2+\j/2) circle (0.75mm);
					\fill (2.5-\i/2,2.5-\j/2) circle (0.75mm);
				}
			}
			
			\draw[-latex] (0,-3) -- (0,3);
			\draw[-latex] (-3,0) -- (3,0);
			\draw[color=red,-latex] (0,0) -- (0.5,0.5);
			\draw[color=red,-latex] (0,0) -- (0.5,1);
			\draw[color=red,dashed] (1.5,3) -- (0,0) -- (3,3);
			\node at (2,1)  {$A^{-1}[Q_{\vec{\ind}}]$};
			\end{tikzpicture}
			\caption{As seen in the figure, we may always assume that $\psi(f)Q_{\vec{\ind}}$ does not contain any quadrant as a strict subset, since, if it does, we may replace $f$ by $f^{-1}$ and obtain a matrix $\psi(f^{-1})$ that does satisfy this condition.}
			\label{fig:replace_f_by_inverse}
		\end{figure}
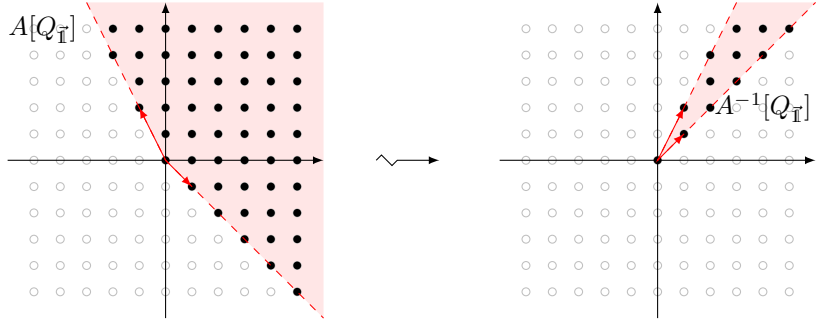
		
		In both cases, if $\psi(f)$ fails to send a quadrant to a quadrant, then either $f$ or $f^{-1}$ cannot be injective, a contradiction since all elements of $\Sym(\shift{\theta},\Z^d)$ are homeomorphisms. Then, since any matrix $\psi(f)$ with $f\in\Sym(\shift{\theta},\Z^d)$ sends quadrants to quadrants, it must send sets that can be written as intersections of quadrants to other sets of the same form. In particular, for any $1\le j\le d$, a set of the form $\pm\N_0\vec{e}_j$ can be written as the intersection of all the quadrants with vertex $\vec{0}$ contained in the half-space containing $\pm\vec{e}_j$; thus, the image of $\pm\N_0\vec{e}_j$ must be another set of this same form. Since $\Z\vec{e}_j=\N_0\vec{e}_j\cup(-\N_0\vec{e}_j)$ and $\psi(f)$ is (a matrix representing) a linear function, this means $\psi(f)$ sends any (discrete) linear subspace $\Z\vec{e}_j$ to some other $\Z\vec{e}_i$, proving the last assertion.\qed
	\end{dem}
	
	The previous results, combined, lead to the following theorem:
	
	\begin{teo}
		\label{teo:structure_of_symmetry_group}
		For a $d$-dimensional, nontrivial, primitive, bijective substitution $\theta$, the quotient group of all admissible lattice transformations of the subshift $\shift{\theta}$,  $\Sym(\shift{\theta},\Z^d)/\Aut(\shift{\theta},\Z^d)$, is isomorphic to a subset of the hyperoctaedral group $Q_d\cong (\Z/2\Z)\wr S_d=(\Z/2\Z)^d\rtimes S_d$, which represents the symmetries of the $d$-dimensional cube. Thus, the extended symmetry group $\Sym(\shift{\theta},\Z^d)$ is virtually-$\Z^d$.
	\end{teo}
	
	\begin{dem}
		As stated previously in Corollary \ref{cor:quadrant_permutations}, since for any $f\in\Sym(\shift{\theta},\Z^d)$ the linear mapping $\psi(f)$ sends quadrants to quadrants, the set $\Sym(\shift{\theta},\Z^d)$ acts on the one-dimensional subspaces $\Z\vec{e}_1,\dots,\Z\vec{e}_d$ by permutation. Thus, $\psi(f)$ must be given by a matrix that sends each $\vec{e}_i$ from the canonical basis to a vector $\pm\vec{e}_j$ and thus each column of $\psi(f)$ is such a vector; this, plus the nonsingularity of $\psi(f)$, shows that the associated matrix must be of the form:
		\[\psi(f)=[(-1)^{t_1}\vec{e}_{\sigma(1)}\mid(-1)^{t_2}\vec{e}_{\sigma(2)}\mid\cdots\mid(-1)^{t_d}\vec{e}_{\sigma(d)} ], \]
		where $\sigma$ is a permutation of $\{1,\dots,d\}$ and $t_1,\dots,t_d\in\{0,1\}$. These matrices correspond to a finite subgroup of $\GL_d(\Z)$ which is isomorphic to $Q_d$. Indeed, the set of all diagonal matrices of this form is isomorphic to $(\Z/2\Z)^d$, while the set of all matrices with nonnegative entries of this form is isomorphic to $S_d$, and any matrix of the aforementioned form is a product of a permutation matrix with positive entries and a diagonal matrix in a unique way.
		
		Thus, $\psi$ can be seen as a group morphism $\Sym(\shift{\theta},\Z^d)\to Q_d$ by identifying the latter with the corresponding matrices; since $\ker(\psi)=\Aut(\shift{\theta},\Z^d)$, we conclude that $\Sym(\shift{\theta},\Z^d)/\Aut(\shift{\theta},\Z^d)\cong\mathrm{im}(\psi)\le Q_d$, the desired result.\qed
	\end{dem}
	
	The previous result imposes a very strict limitation on the structure of the group $\Sym(\shift{\theta},\Z^d)$; thus, with some additional information, we may be able to compute $\Sym(\shift{\theta},\Z^d)$. An example of this is the following:
	
	\begin{cor}
		The extended symmetry group of the generalized Thue-Morse substitution, given by:
		\begin{align*}
			\theta_{\rm TM}:\{0,1\}&\to\{0,1\}^{\{0,1\}^d}\\
			a&\mapsto ((a+m_1+\dots+m_d)\bmod 2)_{(m_1,\dotsc,m_d)\in\{0,1\}^d},
		\end{align*}
		is a semidirect product of the form:
		\[\Sym(\shift{\theta_{TM}},\Z^d)\cong (\Z^d\times\Z/2\Z)\rtimes Q_d, \]
		generated by the shifts, the relabeling map $\delta(x)=\overline{x}$ and the $2^dd!$ rigid symmetries of the coordinate axes given by $(\varphi_A(x))_{\vec{s}}=x_{A\vec{s}}$, with $A\in Q_d$.
	\end{cor}
	
	\begin{dem}
		By Theorem \ref{teo:structure_of_symmetry_group}, $\Sym(\shift{\theta_{\rm TM}},\Z^d)$ is an $\Aut(\shift{\theta_{\rm TM}},\Z^d)$-by-$R$ group extension, where $\Aut(\shift{\theta_{\rm TM}},\Z^d)=\langle\{\sigma_{\vec{k}}\}_{\vec{k}\in\Z^d}\dotcup\{\delta\}\rangle\cong\Z\times(\Z/2\Z)$ and $R$ is a subgroup of $Q_d$. We have to verify, then, that the rigid coordinate symmetries $\varphi_A$ are effectively elements of $\Sym(\shift{\theta_{\rm TM}},\Z^d)$, as they then are mapped to the corresponding elements of $Q_d$ (seen as a matrix group) in $\GL_d(\Z)$ and their composition and inverses are also rigid coordinate symmetries.
		
		Let $f=f_A$ be any rigid coordinate symmetry and let $\Sigma$ be the set of all fixed points of $\theta_{\rm TM}^2$, which are in a $1$-$1$ correspondence with the set of all possible seeds $\{0,1\}^{\{-1,0\}^d}$. By inspection, we see that any rigid coordinate symmetry sends a point of $\Sigma$ to a shift of another point of this set (note that it is possible for the image not to belong to $\Sigma$ as the symbol in the position $\vec{0}$ never changes position); since any subpattern of any point $y\in\shift{\theta}$ is a subpattern of an $x\in\Sigma$, then by the generalized Curtis-Hedlund-Lyndon theorem any subpattern of $f(y)$ is also a subpattern of some other $x'\in\Sigma$; thus, $f(y)\in\shift{\theta}$. We see then that $f$ maps $\shift{\theta}$ to itself; as it has an obvious inverse $f_{A^{-1}}$ which is also a rigid coordinate symmetry, this function $f$ is a homeomorphism $\shift{\theta}\to\shift{\theta}$ satisfying the Curtis-Hedlund-Lyndon condition and hence an element of $\Sym(\shift{\theta},\Z^d)$. Since there is a rigid coordinate symmetry $f_A$ for all $A\in Q_d$, we conclude that the aforementioned group $R$ is all of $Q_d$.
		
		We note that $f_A\circ f_{A'}=f_{AA'}$ and thus $\iota: A\mapsto f_A$ is an embedding of $Q_d$ into $\Sym(\shift{\theta_{\rm TM}},\Z^d)$ and a right inverse for $\psi$. Thus, the short exact sequence $\Aut(\shift{\theta_{\rm TM}},\Z^d)\incl\Sym(\shift{\theta_{\rm TM}},\Z^d)\epi Q_d$ splits, resulting in the desired decomposition of $\Sym(\shift{\theta_{\rm TM}},\Z^d)$ as a semidirect product. \qed
	\end{dem}
	
	\section{The Robinson shift and fractures in subshifts}
	
	In this section, we shall temporarily focus our attention away from the substitutive tilings studied above and analyze a well-known example of \textbf{strongly aperiodic} $\Z^2$-subshift, the Robinson shift.
	
	\begin{defi}
		Let $X$ be a $\Z^d$-subshift. We say $X$ is \textbf{strongly aperiodic} if all points in $X$ have trivial stabilizer, i.e., for all $x\in X$, $\sigma_{\vec{k}}(x)=x$ implies $\vec{k}=\vec{0}$.
	\end{defi}
	
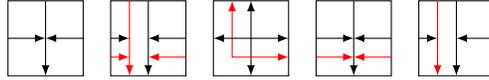
\begin{figure}[h!]
		\centering
		\begin{tikzpicture}
		\draw (0,0) rectangle (1,1);
		\draw[-latex] (0.5,1) -- (0.5,0);
		\draw[-latex] (0,0.5) -- (0.5,0.5);
		\draw[-latex] (1,0.5) -- (0.5,0.5);
		\end{tikzpicture}\quad
		\begin{tikzpicture}
		\draw (0,0) rectangle (1,1);
		\draw[-latex] (0.5,1) -- (0.5,0);
		\draw[-latex,color=red] (0.25,1) -- (0.25,0);
		\draw[-latex] (0,0.5) -- (0.25,0.5);
		\draw[-latex] (1,0.5) -- (0.5,0.5);
		\draw[-latex,color=red] (0,0.25) -- (0.25,0.25);
		\draw[-latex,color=red] (1,0.25) -- (0.5,0.25);
		\end{tikzpicture}\quad
		\begin{tikzpicture}
		\draw (0,0) rectangle (1,1);
		\draw[latex-latex] (0.5,1) -- (0.5,0);
		\draw[latex-latex,color=red] (0.25,1) -- (0.25,0.25) -- (1,0.25);
		\draw[latex-latex] (0,0.5) -- (1,0.5);
		\end{tikzpicture}\quad		
		\begin{tikzpicture}
		\draw (0,0) rectangle (1,1);
		\draw[-latex] (0.5,1) -- (0.5,0);
		\draw[-latex] (0,0.5) -- (0.5,0.5);
		\draw[-latex] (1,0.5) -- (0.5,0.5);
		\draw[-latex,color=red] (0,0.25) -- (0.5,0.25);
		\draw[-latex,color=red] (1,0.25) -- (0.5,0.25);
		\end{tikzpicture}\quad
		\begin{tikzpicture}
		\draw (0,0) rectangle (1,1);
		\draw[-latex] (0.5,1) -- (0.5,0);
		\draw[-latex,color=red] (0.25,1) -- (0.25,0);
		\draw[-latex] (0,0.5) -- (0.25,0.5);
		\draw[-latex] (1,0.5) -- (0.5,0.5);
		\end{tikzpicture}
		\caption{The five types of Robinson tiles, resulting in an alphabet of $28$ symbols after applying all possible rotations and reflections. The third tile is usually called a \textbf{cross}.}
		\label{fig:robinson_tiles}
	\end{figure}

	The Robinson shift is a nearest neighbor two-dimensional shift with added local restrictions (and thus of finite type), whose alphabet consists of all the rotations and reflections of the five tiles from Figure \ref{fig:robinson_tiles}, resulting in $28$ different symbols.
	
	The Robinson shift $X_{\rm Rob}$ is given by the following local rules:
	\begin{enumerate}[label=(\arabic*)]
		\item Every arrow head in a tile must be in contact with an arrow tail from an adjacent tile (nearest neighbor rule). This is similar to the local rule of a Wang tiling (although not exactly equivalent; see \cite{Rob71} or \cite{GJS2012} for details).
		\item There is a translation of the sublattice $2\Z\times2\Z$ that only has rotations of the central tile of Figure \ref{fig:robinson_tiles} (which shall be referred to as \textbf{crosses}).
		\item Any other crosses appear diagonally adjacent to one of the crosses from the sublattice of Rule (2). Namely, if the cross-only sublattice of a given point is $2\Z\times2\Z+\vec{k}$, then any other cross is placed at one of the points from $2\Z\times2\Z+\vec{k}+\vec{\ind}$.
	\end{enumerate}
	It is easy to see that those rules can be enforced with stricty local restrictions and thus $X_{\rm Rob}$ is a shift of finite type. These rules force the $28$ basic tiles to form larger patterns with similar behaviors to each of the five tiles (in particular, patterns of size $(2^n-1)\times(2^n-1)$ that behave as larger analogues of crosses and that are usually referred to as \textbf{$n$-th order supertiles}). By compactness, as we can always build larger supertiles from smaller ones, we can prove that $X_{\rm Rob}$ is a non-empty strongly aperiodic subshift; it is not minimal, but it has a unique minimal subsystem $M_{\rm Rob}$ (which is the factor of a subshift of finite type).
	
	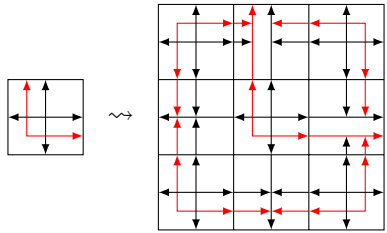
\begin{figure}[!h]
		\centering
		\begin{tikzpicture}
		\draw (0,0) rectangle (1,1);
		\draw[latex-latex] (0.5,1) -- (0.5,0);
		\draw[latex-latex,color=red] (0.25,1) -- (0.25,0.25) -- (1,0.25);
		\draw[latex-latex] (0,0.5) -- (1,0.5);
		
		\draw (2,-1) rectangle (5,2) (2,0) -- (5,0) (2,1) -- (5,1) (3,-1) -- (3,2) (4,-1) -- (4,2);
		\draw[latex-latex,color=red] (2.25,0) -- (2.25,-0.75) -- (3,-0.75);
		\draw[latex-latex,color=red] (2.25,1) -- (2.25,1.75) -- (3,1.75);
		\draw[latex-latex,color=red] (4.75,0) -- (4.75,-0.75) -- (4,-0.75);
		\draw[latex-latex,color=red] (4.75,1) -- (4.75,1.75) -- (4,1.75);
		
		\draw[latex-latex] (2,-0.5) -- (3,-0.5);
		\draw[latex-latex] (2,1.5) -- (3,1.5);
		\draw[latex-latex] (4,-0.5) -- (5,-0.5);
		\draw[latex-latex] (4,1.5) -- (5,1.5);		
		
		\draw[latex-latex] (2.5,-1) -- (2.5,0);
		\draw[latex-latex] (4.5,-1) -- (4.5,0);
		\draw[latex-latex] (2.5,1) -- (2.5,2);
		\draw[latex-latex] (4.5,1) -- (4.5,2);
		
		\draw[latex-latex,color=red] (3.25,1) -- (3.25,0.25) -- (4,0.25);
		\draw[latex-latex] (3,0.5) -- (4,0.5);		
		\draw[latex-latex] (3.5,0) -- (3.5,1);
		
		\draw[-latex] (3.5,1) -- (3.5,2);
		\draw[-latex] (3.5,0) -- (3.5,-1);
		\draw[-latex] (3,0.5) -- (2,0.5);
		\draw[-latex] (4,0.5) -- (5,0.5);
		
		\draw[-latex] (3,1.5) -- (3.25,1.5);
		\draw[-latex] (4,1.5) -- (3.5,1.5);
		\draw[-latex] (3,-0.5) -- (3.5,-0.5);
		\draw[-latex] (4,-0.5) -- (3.5,-0.5);			
		
		\draw[-latex,color=red] (3,1.75) -- (3.25,1.75);
		\draw[-latex,color=red] (4,1.75) -- (3.5,1.75);
		\draw[-latex,color=red] (3,-0.75) -- (3.5,-0.75);
		\draw[-latex,color=red] (4,-0.75) -- (3.5,-0.75);	
		
		\draw[-latex,color=red] (2.25,0) -- (2.25,0.5);
		\draw[-latex,color=red] (2.25,1) -- (2.25,0.5);
		\draw[-latex] (2.5,0) -- (2.5,0.5);
		\draw[-latex] (2.5,1) -- (2.5,0.5);
		
		\draw[-latex,color=red] (4.75,0) -- (4.75,0.25);
		\draw[-latex,color=red] (4.75,1) -- (4.75,0.5);
		\draw[-latex] (4.5,0) -- (4.5,0.25);
		\draw[-latex] (4.5,1) -- (4.5,0.5);
		
		\draw[-latex,color=red] (3.25,1) -- (3.25,2);
		\draw[-latex,color=red] (4,0.25) -- (5,0.25);
		
		\node at (1.5,0.5) {$\rightsquigarrow$};			
		\end{tikzpicture}
		\caption{The formation of a second order supertile of size $3\times 3$.}
		\label{fig:supertile}
	\end{figure}
	
	The following result has been proven by Sebastián Donoso and Wenbo Sun, in \cite{DS2014}:
	
	\begin{teo}
		$\Aut(M_{\rm Rob},\Z^2)=\langle\sigma_{(1,0)},\sigma_{(0,1)}\rangle\cong\Z^2$.
	\end{teo}
	
	From this result it is possible to show that the same holds for $X_{\rm Rob}$, namely that the only automorphisms of the Robinson shift are the trivial ones. We aim to extend this result by computing the extended symmetry group of the Robinson shift. For this, we need to introduce a distinguished subset of $\Z^2$ which represents part of the structure of a shift which is preserved by extended symmetries:
	
	\begin{defi}
		Let $X$ be a strongly aperiodic $\Z^2$-subshift. We say $X$ has a \textbf{fracture} in the direction $\vec{q}\in\Z^2$ if there is a point $x^*\in X$, infinite different values $k_1<k_2<k_3<\dotsc\in\Z$, and two disjoint half-planes $S^+,S^-\subseteq\Z^2$ separated by $\Z\vec{q}$ (i.e. $S^+\cap S^-=S^+\cap\Z\vec{q}=S^-\cap\Z\vec{q}=\varnothing$; it is not necessary that $S^+\cup S^-\cup\Z\vec{q}=\Z^2$) such that, for each $j\in\N$, there is a point $x^{(j)}\in X$ that satisfies the two conditions:
		\[x^{(j)}|_{S^+}=x^*|_{S^+},\quad x^{(j)}|_{S^-}=\sigma_{k_j\vec{q}}(x^*)|_{S^-}. \]
	\end{defi}
	
	\begin{nota}
		We exclude subshifts with periodic points from this definition as, if $x\in\Per_{\vec{p}}(X)$, then we may take $k_j=j$ and $x^{(j)}=x$ for all values of $j$, resulting in a point with a fracture in the direction $\vec{p}$; this makes the definition of direction of fracture redundant with the concept of direction of periodicity, which is also preserved by extended symmetries.
	\end{nota}
	
	\begin{lem}
		Let $\vec{q}$ be a direction of fracture for a two-dimensional strongly aperiodic subshift $X$ and $f\in\Sym(X,\Z^2)$. Then $\psi(f)\vec{q}$ is a direction of fracture as well. 
	\end{lem}
	
	\begin{dem}
		Let $\vec{q}$ be a direction of fracture, and $x^*, (x^{(j)})_{j\in\N},(k_j)_{j\in\N}$ be the associated points and magnitudes from the definition above. By the generalized Curtis-Hedlund-Lyndon theorem, as $x^*|_{S^+}=x^{(j)}|_{S^+}$, then $f(x^*)|_{\psi(f)((S^+)^{\circ r})}=f(x^{(j)})|_{\psi(f)((S^+)^{\circ r})}$, where $r$ is the radius of the symmetry $f$. By the same argument, and since $f\circ\sigma_{\vec{q}}=\sigma_{\psi(f)\vec{q}}\circ f$, we conclude that $f(x^{(j)})|_{\psi(f)((S^-)^{\circ r})}=\sigma_{k_j\psi(f)\vec{q}}\circ f(x^*)|_{\psi(f)((S^-)^{\circ r})}$.
		
		Note that, since $S^+$ and $S^-$ are half-planes disjoint from the linear subspace $\Z\vec{q}$, and $\psi(f)$ is a linear map, $(S^\pm)^{\circ r}$ are also half-spaces and thus their corresponding images $\psi(f)((S^\pm)^\circ r)$ are half-spaces as well. As subsets of the images of disjoint sets, they are also disjoint from $\Z(\psi(f)\vec{q})$ and from each other. Thus, by defining $y^*=f(x^*),y^{(j)}=f(x^{(j)})$ we see that these points conform a fracture of $X$ in the direction $\psi(f)\vec{q}$.\qed
	\end{dem}
	
	This result provides a subset of $\Z^2$ over which the set $\Sym(X,\Z^2)$ is forced to act ``naturally''; thus, if this subset has sufficiently strong constraints coming from the structure of $X$, this enforces similar restrictions on the possible values of $\psi(f)$ for $f\in\Sym(X,\Z^2)$. As we shall see below, this is the case for the aperiodic Robinson shift:
	
	\begin{prop}
		For the Robinson shift, $\Sym(X_{\rm Rob},\Z^2)\cong\Z^2\rtimes D_4$, where $D_4$ is the dihedral group of order $8$.
	\end{prop}
	
	\begin{dem}
		To prove this result, we will show that the set $\mathcal{S}$ of all directions of fracture of $X_{\rm Rob}$ is $\Z\vec{e}_1\cup\Z\vec{e}_2$. Assuming this as true, we see that, since $\psi(f)$ is always a $\Z$-invertible matrix, it must send $\{\vec{e}_1,\vec{e}_2\}$ to a basis of $\Z^2$ contained in $\Z\vec{e}_1\cup\Z\vec{e}_2$, which is always a two-element set of the form $\{\pm\vec{e}_1,\pm\vec{e}_2 \}$ or $\{\pm\vec{e}_1,\mp\vec{e}_2\}$, and thus the elements of $\Sym(X_{\rm Rob},\Z^2)$ correspond to one of the eight possible matrices belonging to the standard copy of $D_4=Q_2$ (defined in the previous section) in $\GL_2(\Z)$. Then, by finding an explicit subgroup of $\Sym(X_{\rm Rob},\Z^2)$ isomorphic to $D_4$ by $\psi$, we deduce the claimed semidirect product decomposition.
		
		To show that $X_{\rm Rob}$ has fractures in the directions $\vec{e}_1$ and $\vec{e}_2$, we need to recall some basic details about the construction of an infinite valid configuration of the Robinson shift. As stated above, the five basic Robinson tiles (together with their rotations and reflections) combine to form $3\times 3$ patterns with a similar behavior to crosses, named second order supertiles. Four of these second order supertiles, together with smaller substructures, further combine to form $7\times 7$ patterns (third order supertiles) and so on. In every case, the central tile of an $n$-th order supertile is a cross, which gives an orientation to the supertile in a similar way to the two-headed, L-shaped arrow on a cross.
		
		We may fill the whole upper right quadrant $Q_{\vec{\ind}}=\N^2$ as follows: we start by placing a cross on its vertex $\vec{0}$ with its L-shaped arrow pointing up and right, and then place another cross with the same orientation at the position $(1,1)$. This new cross, together with the previously placed one, allows us to fill the lower $3\times 3$ section of $\N^2$, $[0,2]^2$, with a second order supertile. We iterate this process by placing a cross with the same orientation at the position $(3,3),(7,7),\dots,(2^n-1,2^n-1),\dots$ and constructing the corresponding second, third, \dots $n$-th order supertile and so on. By compactness, there is only one way to fill all of $\N^2$ as a limit to this process. We call the configuration obtained an infinite order supertile.
		
		We may fill the other three quadrants with similar constructions resulting in infinite order supertiles with different orientations, each of these separated from the other infinite supertiles by a row or column of copies of the first tile from Figure \ref{fig:robinson_tiles}. As we see in Figure \ref{fig:robinson_supertile_separation}, this will result in a translate of $\Z\times\{0\}\cup\{0\}\times\Z$ containing only copies of this tile, with all of the tiles in one of the strips $\Z\times\{0\}$ or $\{0\}\times\Z$ (the latter in the figure) having the same orientation, while the other strip will have all of its tiles pointing towards the center.
		
		\begin{figure}[!ht]
			\centering
			\includegraphics[scale=0.6]{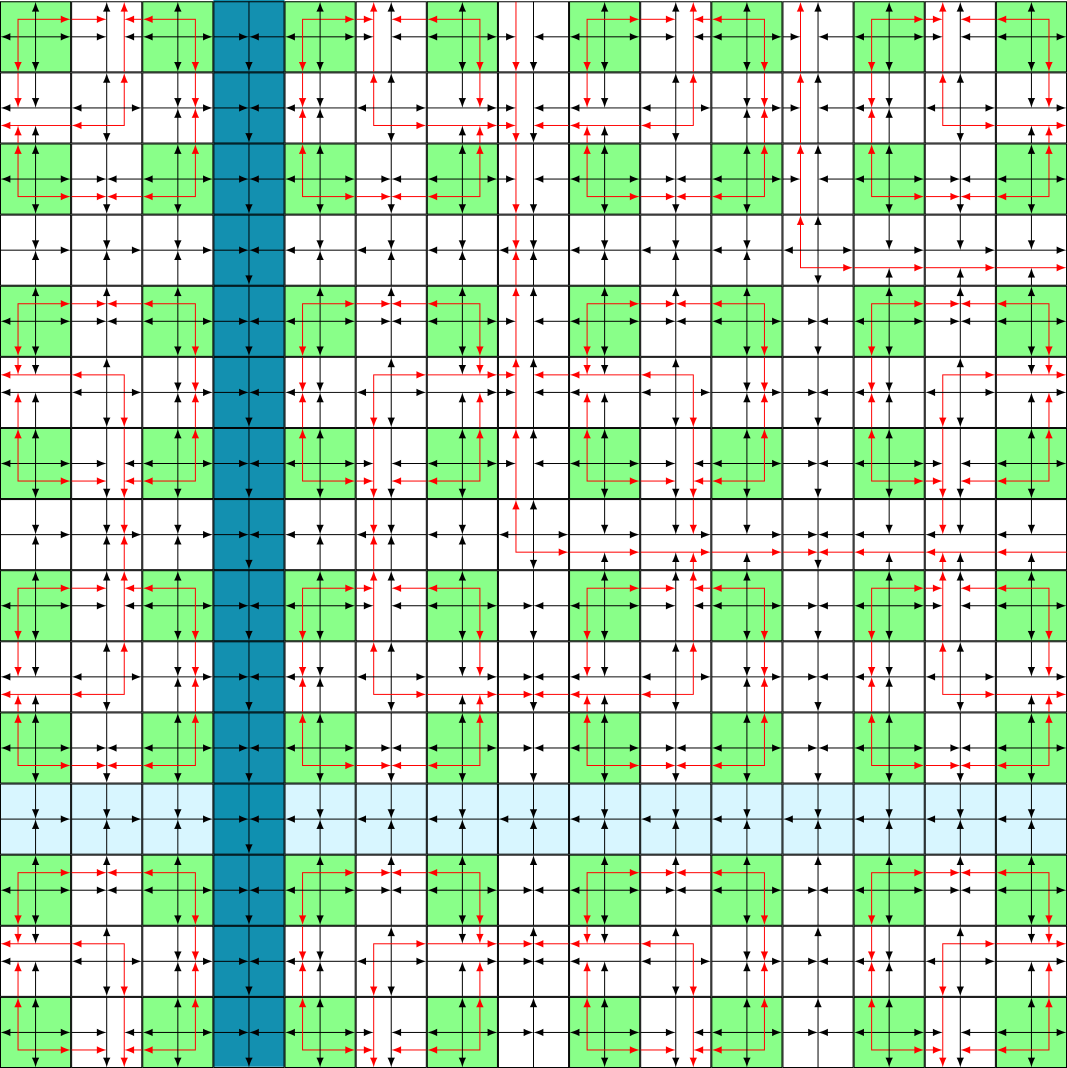}
			\caption{A fragment of a point from the Robinson shift, distinguishing the four supertiles involved, the vertical and horizontal strips of tiles separating each supertile and the $2\Z\times 2\Z$ sublattice that contains only crosses. Note that the tiles in the vertical strip separating the supertiles are copies of the first tile of Figure \ref{fig:robinson_tiles} with the same orientation.}
			\label{fig:robinson_supertile_separation}
		\end{figure}
		
		Since the Robinson shift behaves like a nearest-neighbor shift with added restrictions, the existence of a vertical (resp., horizontal) strip with copies of the same tile allows us to vertically shift the tiles contained in the right half-plane however we see fit, as long as the coset of the $2\Z\times 2\Z$ sublattice containing only crosses is respected. In practice, this shows that in the point $x\in X_{\rm Rob}$ represented partially in Figure \ref{fig:robinson_supertile_separation} we may replace the tiles from the right half-plane with the corresponding tiles from $\sigma_{(0,2k)}(x)$ and obtain valid points for all values of $k\in\Z$. We see an example of this in Figure \ref{fig:robinson_breaking}.
		\begin{figure}[!ht]
			\centering
			\includegraphics[scale=0.6]{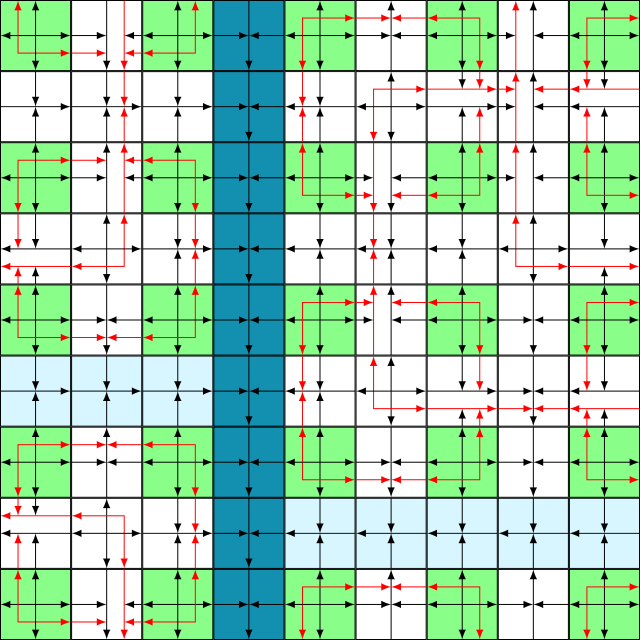}\qquad
			\includegraphics[scale=0.6]{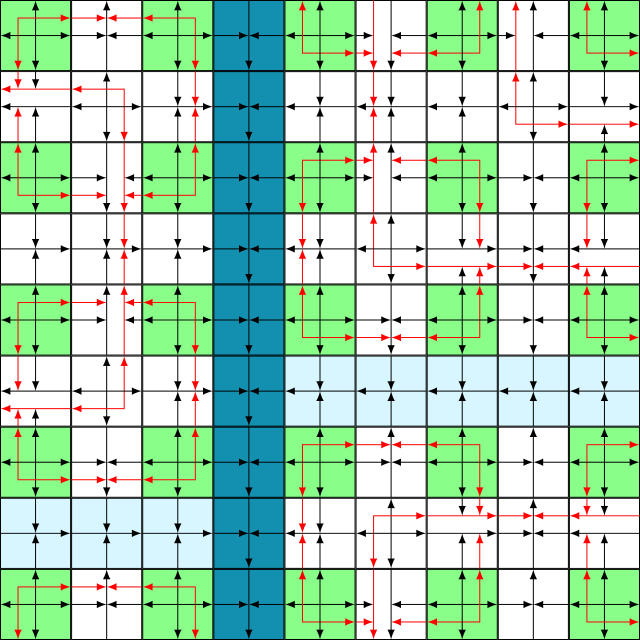}
			\caption{Two possible ways in which the tiling from Figure \ref{fig:robinson_supertile_separation} exhibits fracture-like behavior while still resulting in a valid point from $X_{\rm Rob}$.}
			\label{fig:robinson_breaking}
		\end{figure}
		
		This procedure shows that $X_{\rm Rob}$ has $\vec{e}_1$ and $\vec{e}_2$ as directions of fracture. Now, we need to show that all directions of fracture are contained in the set $\Z\vec{e}_1\cup\Z\vec{e}_2$, and thus all matrices from $\psi[\Sym(X_{\rm Rob},\Z^2)]$ send the set $\{\vec{e}_1,\vec{e}_2\}$ to a linearly independent subset of $\{\vec{e}_1,\vec{e}_2,-\vec{e}_1,-\vec{e}_2\}$. The argument we shall use for this follows a similar outline to the technique used in the proof of Corollary \ref{cor:quadrant_permutations}: the points of the Robinson shift form a hierarchical structure away from a horizontal or vertical fracture, allowing for a decomposition into subpatterns of arbitrarily large size $s$	placed correlative to a lattice of the form $2^n\Z\times 2^n\Z$ (this is similar to the decomposition of a point from a substitutive subshift into patterns of the form $\theta^m(a),a\in\mathcal{A}$ for arbitrarily large values of $m$). The existence of fractures that are neither vertical nor horizontal would result in ``ruptures'' in this hierarchical structure, leading to a contradiction.
		
		Formally, we proceed as follows. Suppose that $X_{\rm Rob}$ has a fracture in the direction $\vec{q}\in\Z^2\setminus(\Z\vec{e}_1\cup\Z\vec{e}_2)$, and let $S^+,S^-$ be the disjoint half-planes separated by $\vec{q}$. The set $F_{\vec{q}}=\Z^2\setminus (S^+\dotcup S^-)$ is necessarily of the form $\Z\vec{q}+[\vec{r}_1,\vec{r}_2]$, namely, a finite union of translates of $\Z\vec{q}$, and thus its intersection with any set of the form $\Z\times\{k\}$ or $\{k\}\times\Z$ is finite. This is because the intersection of such a set with $\Z\vec{q}$ consists of at most a single point, as $\vec{q}$ is not a multiple of $\vec{e}_1$ nor $\vec{e}_2$. Thus, for any sufficiently large value $M\in\N$, it is easy to verify that for any point $\vec{p}\in F_{\vec{q}}$, any translation of the rectangle $[-M\vec{\ind},M\vec{\ind}]$ that contains $\vec{p}$ also contains points from either $S^+$ or $S^-$ (or both).
		
		\begin{figure}[!ht]
			\centering
			\begin{tikzpicture}[scale=0.6]
			\fill[color=red!20!white] (-4,-4) -- (-4,-3) -- (3,4) -- (4,4) -- (4,3) -- (-3,-4) -- cycle;
			\draw[-latex] (-4,0) -- (4,0);
			\draw[-latex] (0,-4) -- (0,4);
			\draw[color=red] (-4,-4) -- (4,4);
			\draw[color=red,dashed] (-4,-3) -- (3,4) (-3,-4) -- (4,3);
			\foreach \i in {-4,-2,0,2} {
				\foreach \j in {-4,-2,0,2} {
					\draw[color=blue] (\i+0.1,\j+0.1) rectangle (\i+1.9,\j+1.9);
				}
			}
			\node at (-1,3) {$S^+$};
			\node at (3,-1) {$S^-$};
			\node at (4,4) [right] {$\Z\vec{q}$};
			\node at (-4,-4) [left,color=red] {$F_{\vec{q}}$};
			\end{tikzpicture}
			\caption{The substructure of a point of $X_{\rm Rob}$ in terms of $n$-th order supertiles. Note how all supertiles overlap either $S^+$ or $S^-$.}
		\end{figure}
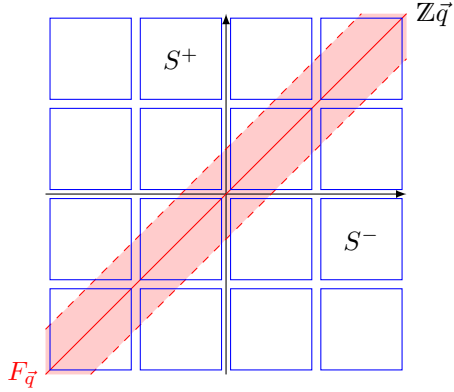
		
		Choose $n\in\N,n>1$ such that for $M=2^n-1$ the above condition holds, while satisfying the additional condition $M>2k_1\|\vec{q}\|_1$. All $n$-th order supertiles, thus, contain points from $S^+\dotcup S^-$. Let $\vec{p}$ be an element from $F_{\vec{q}}$ that belongs to the support of an $n$-th order supertile, and suppose this supertile overlaps the half-plane $S^+$. If there is no such supertile, all $n$-th order supertiles containing points of $F_{\vec{q}}$ only overlap $S^-$, implying, since $S^+$ is the intersection of a real half-plane $H_{\vec{\alpha},c}=\{\vec{v}\in\R^2:\langle\vec{v},\vec{\alpha}\rangle \ge c \}$ with $\Z^2$, that $S^+$ is a translation of $\Z\times(\pm\N)$ (or $(\pm\N)\times\Z$). C
		onvex combinations of points of $S^+$ with integer coefficients belong to $S^+$ as well, so $S^+$ cannot have ``gaps'', and it is a union of disjoint, horizontally or vertically adjacent translates of $[1,2^n]^2$. This implies that $\vec{q}$ is in the set $\Z\vec{e}_1\cup\Z\vec{e}_2$, a contradiction; thus, the aforementioned supertile exists. Evidently, the same argument shows the existence of other $n$-th order supertiles which intersect $S^{-}$.
		
		Since each horizontal or vertical strip $F_{\vec{q}}\cap(\Z\times\{k\})$ (resp. $F_{\vec{q}}\cap(\Z\times\{k\})$) intersects finitely many supertiles, we see that the arrangement of the $n$-th order supertiles in $S^+$ away from a vertical or horizontal fracture (which in this case must correspond to a bi-infinite column or row of copies of the first tile from Figure \ref{fig:robinson_tiles}, all with the same orientation) affects the placement of the supertiles in $S^-$ as well. However, since the tiling has a fracture in the direction $\vec{q}$, we may shift the supertiles in $S^-$ by $k_1\vec{q}$ and obtain a valid configuration. By our choice of $n$, the shift $\sigma_{k_1\vec{q}}$ moves the $n$-th order supertiles by less than $M$ units both horizontally and vertically (since $M>2k_1\|\vec{q}\|_1$), and thus the supertiles in $S^-$ are shifted to a position that does not match the arrangement of supertiles from $S^+$. We may see this situation in Figure \ref{fig:supertile_breaker}.
		
		\begin{figure}[!h]
			\centering
			\begin{tikzpicture}[scale=0.6]
			\fill[color=red!20!white] (-4,-4) -- (-4,-3) -- (3,4) -- (4,4) -- (4,3) -- (-3,-4) -- cycle;
			\draw[-latex] (-4,0) -- (4,0);
			\draw[-latex] (0,-4) -- (0,4);
			\draw[color=red] (-4,-4) -- (4,4);
			\draw[color=red,dashed] (-4,-3) -- (3,4) (-3,-4) -- (4,3);
			\foreach \i in {-4,-2,0,2} {
				\foreach \j in {-4,-2,0,2} {
					\draw[color=blue!40!white] (\i+0.1,\j+0.1) rectangle (\i+1.9,\j+1.9);
				}
			}
			\draw[color=green!60!black] (-2.5,-3.5) -- (-1.7,-3.5) -- (-1.7,-2.7);
			\draw[color=green!60!black] (-1.5,-2.5) -- (-1.5,-3.5) -- (0.3,-3.5) -- (0.3,-1.7) -- (-0.7,-1.7);
			\draw[color=green!60!black] (-0.5,-1.5) -- (0.3,-1.5) -- (0.3,-0.7);
			\draw[color=green!60!black] (0.5,-0.5) -- (0.5,-1.5) -- (2.3,-1.5) -- (2.3,0.3) -- (1.3,0.3);
			\draw[color=green!60!black] (1.5,0.5) -- (2.3,0.5) -- (2.3,1.3);
			\draw[color=green!60!black] (2.5,1.5) -- (2.5,0.5) -- (4,0.5) (3.3,2.3) -- (4,2.3) (3.5,2.5) -- (4,2.5);
			\draw[color=green!60!black] (0.5,-3.5) rectangle +(1.8,1.8);
			\draw[color=green!60!black] (4,0.3) -- (2.5,0.3) -- (2.5,-1.5) -- (4,-1.5);
			\draw[color=green!60!black] (4,-1.7) -- (2.5,-1.7) -- (2.5,-3.5) -- (4,-3.5);
			\draw[color=green!60!black] (-2.7,-3.7) -- (-1.7,-3.7) -- (-1.7,-4) (-1.5,-4) -- (-1.5,-3.7) -- (0.3,-3.7) -- (0.3,-4) (0.5,-4) -- (0.5,-3.7) -- (2.3,-3.7) -- (2.3,-4) (2.5,-4) -- (2.5,-3.7) -- (4,-3.7);
			
			\node at (2,-2) [rotate=45,color=red,above] {$\sigma_{k_1\vec{q}}$};
			\draw[color=red,-latex] (1,-3) -- (3,-1);
			\node at (-3,3) {$S^+$};
			\node at (3,-3) {$S^-$};
			\node at (4,4) [right] {$\Z\vec{q}$};
			\node at (-4,-4) [left,color=red] {$F_{\vec{q}}$};
			\end{tikzpicture}
			\caption{How a shift by $k_1\vec{q}$ makes the arrangement of supertiles in $S^+$ not match with the corresponding tiles in $S^-$.}
			\label{fig:supertile_breaker}
		\end{figure}
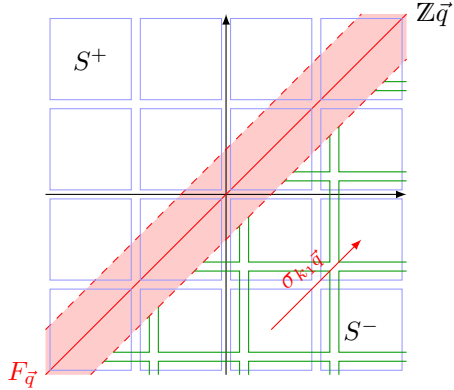
		
		Given that we are assuming that this point (say, $x$) is a fracture point for $X_{\rm Rob}$, there must be some other point $y$ which matches $x$ in $S^+$ and $\sigma_{k_1\vec{q}}(x)$ in $S^-$, which breaks the rigidity of the structure of supertiles imposed by the rules of the Robinson shift. Thus, fractures along non-principal directions cannot exist.
		
		Finally, we need to construct a copy of $D_4$ contained in $\Sym(X_{\rm Rob},\Z^2)$. For this, since $D_4$ is a $2$-generated group, we only need to show the existence of two extended symmetries $\rho,\mu:X_{\rm Rob}\to X_{\rm Rob}$, mapped respectively by $\psi$ to the matrices:
		\[\psi(\rho)=\begin{bmatrix}
		0 & -1 \\
		1 & 0
		\end{bmatrix},\quad \psi(\mu)=\begin{bmatrix}
		-1 & 0 \\
		0 & 1
		\end{bmatrix}, \]
		since these two matrices generate an isomorphic copy of $D_4$ contained in $\GL_2(\Z)$. These symmetries $\rho$ and $\mu$ are essentially rigid symmetries of the coordinate axes; however, a composition with a relabeling map is also needed, to replace every tile with the corresponding reflection or rotation. For instance, if we define $\mathfrak{R}:\mathcal{A}\to\mathcal{A}$ as the mapping which assigns to each of the $28$ symbols its corresponding rotation by $\frac{1}{2}\pi$, as seen in Figure \ref{fig:rotation_tiles}, then $\rho(x)_{(i,j)}=\mathfrak{R}(x_{(-j,i)})$ is the desired symmetry. In the same way, by defining $\mathfrak{M}:\mathcal{A}\to\mathcal{A}$ as the mapping that sends each tile to its reflection through the horizontal axes, then we define $\mu$ by the relation $\mu(x)_{(i,j)}=\mathfrak{M}(x_{(-i,j)})$.
		
		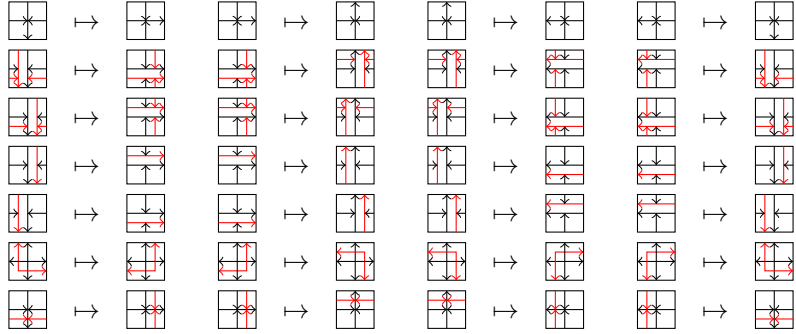
\begin{figure}[!ht]
			\centering
			\[\begin{array}{rclcrclcrclcrcl}
			\begin{tikzpicture}[scale=0.5]
			\draw (0,0) rectangle (1,1);
			\draw[->] (0.5,1) -- (0.5,0);
			\draw[->] (0,0.5) -- (0.5,0.5);
			\draw[->] (1,0.5) -- (0.5,0.5);
			\end{tikzpicture} &\raisebox{1ex}{$\mapsto$} &\begin{tikzpicture}[scale=0.5,rotate=90]
			\draw (0,0) rectangle (1,1);
			\draw[->] (0.5,1) -- (0.5,0);
			\draw[->] (0,0.5) -- (0.5,0.5);
			\draw[->] (1,0.5) -- (0.5,0.5);
			\end{tikzpicture} & &
			\begin{tikzpicture}[scale=0.5,rotate=90]
			\draw (0,0) rectangle (1,1);
			\draw[->] (0.5,1) -- (0.5,0);
			\draw[->] (0,0.5) -- (0.5,0.5);
			\draw[->] (1,0.5) -- (0.5,0.5);
			\end{tikzpicture} &\raisebox{1ex}{$\mapsto$} &\begin{tikzpicture}[scale=0.5,rotate=180]
			\draw (0,0) rectangle (1,1);
			\draw[->] (0.5,1) -- (0.5,0);
			\draw[->] (0,0.5) -- (0.5,0.5);
			\draw[->] (1,0.5) -- (0.5,0.5);
			\end{tikzpicture} & &
			\begin{tikzpicture}[scale=0.5,rotate=180]
			\draw (0,0) rectangle (1,1);
			\draw[->] (0.5,1) -- (0.5,0);
			\draw[->] (0,0.5) -- (0.5,0.5);
			\draw[->] (1,0.5) -- (0.5,0.5);
			\end{tikzpicture} & \raisebox{1ex}{$\mapsto$} & \begin{tikzpicture}[scale=0.5,rotate=270]
			\draw (0,0) rectangle (1,1);
			\draw[->] (0.5,1) -- (0.5,0);
			\draw[->] (0,0.5) -- (0.5,0.5);
			\draw[->] (1,0.5) -- (0.5,0.5);
			\end{tikzpicture} & &
			\begin{tikzpicture}[scale=0.5,rotate=270]
			\draw (0,0) rectangle (1,1);
			\draw[->] (0.5,1) -- (0.5,0);
			\draw[->] (0,0.5) -- (0.5,0.5);
			\draw[->] (1,0.5) -- (0.5,0.5);
			\end{tikzpicture} &\raisebox{1ex}{$\mapsto$} &\begin{tikzpicture}[scale=0.5]
			\draw (0,0) rectangle (1,1);
			\draw[->] (0.5,1) -- (0.5,0);
			\draw[->] (0,0.5) -- (0.5,0.5);
			\draw[->] (1,0.5) -- (0.5,0.5);
			\end{tikzpicture} \\
			\begin{tikzpicture}[scale=0.5]
			\draw (0,0) rectangle (1,1);
			\draw[->] (0.5,1) -- (0.5,0);
			\draw[->,color=red] (0.25,1) -- (0.25,0);
			\draw[->] (0,0.5) -- (0.25,0.5);
			\draw[->] (1,0.5) -- (0.5,0.5);
			\draw[->,color=red] (0,0.25) -- (0.25,0.25);
			\draw[->,color=red] (1,0.25) -- (0.5,0.25);
			\end{tikzpicture}& \raisebox{1ex}{$\mapsto$} &
			\begin{tikzpicture}[scale=0.5,rotate=90]
			\draw (0,0) rectangle (1,1);
			\draw[->] (0.5,1) -- (0.5,0);
			\draw[->,color=red] (0.25,1) -- (0.25,0);
			\draw[->] (0,0.5) -- (0.25,0.5);
			\draw[->] (1,0.5) -- (0.5,0.5);
			\draw[->,color=red] (0,0.25) -- (0.25,0.25);
			\draw[->,color=red] (1,0.25) -- (0.5,0.25);
			\end{tikzpicture} & & 
			\begin{tikzpicture}[scale=0.5,rotate=90]
			\draw (0,0) rectangle (1,1);
			\draw[->] (0.5,1) -- (0.5,0);
			\draw[->,color=red] (0.25,1) -- (0.25,0);
			\draw[->] (0,0.5) -- (0.25,0.5);
			\draw[->] (1,0.5) -- (0.5,0.5);
			\draw[->,color=red] (0,0.25) -- (0.25,0.25);
			\draw[->,color=red] (1,0.25) -- (0.5,0.25);
			\end{tikzpicture}& \raisebox{1ex}{$\mapsto$} &
			\begin{tikzpicture}[scale=0.5,rotate=180]
			\draw (0,0) rectangle (1,1);
			\draw[->] (0.5,1) -- (0.5,0);
			\draw[->,color=red] (0.25,1) -- (0.25,0);
			\draw[->] (0,0.5) -- (0.25,0.5);
			\draw[->] (1,0.5) -- (0.5,0.5);
			\draw[->,color=red] (0,0.25) -- (0.25,0.25);
			\draw[->,color=red] (1,0.25) -- (0.5,0.25);
			\end{tikzpicture} & &
			\begin{tikzpicture}[scale=0.5,rotate=180]
			\draw (0,0) rectangle (1,1);
			\draw[->] (0.5,1) -- (0.5,0);
			\draw[->,color=red] (0.25,1) -- (0.25,0);
			\draw[->] (0,0.5) -- (0.25,0.5);
			\draw[->] (1,0.5) -- (0.5,0.5);
			\draw[->,color=red] (0,0.25) -- (0.25,0.25);
			\draw[->,color=red] (1,0.25) -- (0.5,0.25);
			\end{tikzpicture}& \raisebox{1ex}{$\mapsto$} &
			\begin{tikzpicture}[scale=0.5,rotate=270]
			\draw (0,0) rectangle (1,1);
			\draw[->] (0.5,1) -- (0.5,0);
			\draw[->,color=red] (0.25,1) -- (0.25,0);
			\draw[->] (0,0.5) -- (0.25,0.5);
			\draw[->] (1,0.5) -- (0.5,0.5);
			\draw[->,color=red] (0,0.25) -- (0.25,0.25);
			\draw[->,color=red] (1,0.25) -- (0.5,0.25);
			\end{tikzpicture} & &
			\begin{tikzpicture}[scale=0.5,rotate=270]
			\draw (0,0) rectangle (1,1);
			\draw[->] (0.5,1) -- (0.5,0);
			\draw[->,color=red] (0.25,1) -- (0.25,0);
			\draw[->] (0,0.5) -- (0.25,0.5);
			\draw[->] (1,0.5) -- (0.5,0.5);
			\draw[->,color=red] (0,0.25) -- (0.25,0.25);
			\draw[->,color=red] (1,0.25) -- (0.5,0.25);
			\end{tikzpicture}& \raisebox{1ex}{$\mapsto$} &
			\begin{tikzpicture}[scale=0.5]
			\draw (0,0) rectangle (1,1);
			\draw[->] (0.5,1) -- (0.5,0);
			\draw[->,color=red] (0.25,1) -- (0.25,0);
			\draw[->] (0,0.5) -- (0.25,0.5);
			\draw[->] (1,0.5) -- (0.5,0.5);
			\draw[->,color=red] (0,0.25) -- (0.25,0.25);
			\draw[->,color=red] (1,0.25) -- (0.5,0.25);
			\end{tikzpicture} \\
			\begin{tikzpicture}[xscale=-0.5,yscale=0.5]
			\draw (0,0) rectangle (1,1);
			\draw[->] (0.5,1) -- (0.5,0);
			\draw[->,color=red] (0.25,1) -- (0.25,0);
			\draw[->] (0,0.5) -- (0.25,0.5);
			\draw[->] (1,0.5) -- (0.5,0.5);
			\draw[->,color=red] (0,0.25) -- (0.25,0.25);
			\draw[->,color=red] (1,0.25) -- (0.5,0.25);
			\end{tikzpicture}& \raisebox{1ex}{$\mapsto$} &
			\begin{tikzpicture}[xscale=-0.5,yscale=0.5,rotate=-90]
			\draw (0,0) rectangle (1,1);
			\draw[->] (0.5,1) -- (0.5,0);
			\draw[->,color=red] (0.25,1) -- (0.25,0);
			\draw[->] (0,0.5) -- (0.25,0.5);
			\draw[->] (1,0.5) -- (0.5,0.5);
			\draw[->,color=red] (0,0.25) -- (0.25,0.25);
			\draw[->,color=red] (1,0.25) -- (0.5,0.25);
			\end{tikzpicture} & & 
			\begin{tikzpicture}[xscale=-0.5,yscale=0.5,rotate=-90]
			\draw (0,0) rectangle (1,1);
			\draw[->] (0.5,1) -- (0.5,0);
			\draw[->,color=red] (0.25,1) -- (0.25,0);
			\draw[->] (0,0.5) -- (0.25,0.5);
			\draw[->] (1,0.5) -- (0.5,0.5);
			\draw[->,color=red] (0,0.25) -- (0.25,0.25);
			\draw[->,color=red] (1,0.25) -- (0.5,0.25);
			\end{tikzpicture}& \raisebox{1ex}{$\mapsto$} &
			\begin{tikzpicture}[xscale=-0.5,yscale=0.5,rotate=180]
			\draw (0,0) rectangle (1,1);
			\draw[->] (0.5,1) -- (0.5,0);
			\draw[->,color=red] (0.25,1) -- (0.25,0);
			\draw[->] (0,0.5) -- (0.25,0.5);
			\draw[->] (1,0.5) -- (0.5,0.5);
			\draw[->,color=red] (0,0.25) -- (0.25,0.25);
			\draw[->,color=red] (1,0.25) -- (0.5,0.25);
			\end{tikzpicture} & &
			\begin{tikzpicture}[xscale=-0.5,yscale=0.5,rotate=180]
			\draw (0,0) rectangle (1,1);
			\draw[->] (0.5,1) -- (0.5,0);
			\draw[->,color=red] (0.25,1) -- (0.25,0);
			\draw[->] (0,0.5) -- (0.25,0.5);
			\draw[->] (1,0.5) -- (0.5,0.5);
			\draw[->,color=red] (0,0.25) -- (0.25,0.25);
			\draw[->,color=red] (1,0.25) -- (0.5,0.25);
			\end{tikzpicture}& \raisebox{1ex}{$\mapsto$} &
			\begin{tikzpicture}[xscale=-0.5,yscale=0.5,rotate=-270]
			\draw (0,0) rectangle (1,1);
			\draw[->] (0.5,1) -- (0.5,0);
			\draw[->,color=red] (0.25,1) -- (0.25,0);
			\draw[->] (0,0.5) -- (0.25,0.5);
			\draw[->] (1,0.5) -- (0.5,0.5);
			\draw[->,color=red] (0,0.25) -- (0.25,0.25);
			\draw[->,color=red] (1,0.25) -- (0.5,0.25);
			\end{tikzpicture} & &
			\begin{tikzpicture}[xscale=-0.5,yscale=0.5,rotate=-270]
			\draw (0,0) rectangle (1,1);
			\draw[->] (0.5,1) -- (0.5,0);
			\draw[->,color=red] (0.25,1) -- (0.25,0);
			\draw[->] (0,0.5) -- (0.25,0.5);
			\draw[->] (1,0.5) -- (0.5,0.5);
			\draw[->,color=red] (0,0.25) -- (0.25,0.25);
			\draw[->,color=red] (1,0.25) -- (0.5,0.25);
			\end{tikzpicture}& \raisebox{1ex}{$\mapsto$} &
			\begin{tikzpicture}[xscale=-0.5,yscale=0.5]
			\draw (0,0) rectangle (1,1);
			\draw[->] (0.5,1) -- (0.5,0);
			\draw[->,color=red] (0.25,1) -- (0.25,0);
			\draw[->] (0,0.5) -- (0.25,0.5);
			\draw[->] (1,0.5) -- (0.5,0.5);
			\draw[->,color=red] (0,0.25) -- (0.25,0.25);
			\draw[->,color=red] (1,0.25) -- (0.5,0.25);
			\end{tikzpicture}\\
			\begin{tikzpicture}[xscale=-0.5,yscale=0.5]
			\draw (0,0) rectangle (1,1);
			\draw[->] (0.5,1) -- (0.5,0);
			\draw[->,color=red] (0.25,1) -- (0.25,0);
			\draw[->] (0,0.5) -- (0.25,0.5);
			\draw[->] (1,0.5) -- (0.5,0.5);
			\end{tikzpicture}& \raisebox{1ex}{$\mapsto$} &
			\begin{tikzpicture}[xscale=-0.5,yscale=0.5,rotate=-90]
			\draw (0,0) rectangle (1,1);
			\draw[->] (0.5,1) -- (0.5,0);
			\draw[->,color=red] (0.25,1) -- (0.25,0);
			\draw[->] (0,0.5) -- (0.25,0.5);
			\draw[->] (1,0.5) -- (0.5,0.5);
			\end{tikzpicture} & & 
			\begin{tikzpicture}[xscale=-0.5,yscale=0.5,rotate=-90]
			\draw (0,0) rectangle (1,1);
			\draw[->] (0.5,1) -- (0.5,0);
			\draw[->,color=red] (0.25,1) -- (0.25,0);
			\draw[->] (0,0.5) -- (0.25,0.5);
			\draw[->] (1,0.5) -- (0.5,0.5);
			\end{tikzpicture}& \raisebox{1ex}{$\mapsto$} &
			\begin{tikzpicture}[xscale=-0.5,yscale=0.5,rotate=180]
			\draw (0,0) rectangle (1,1);
			\draw[->] (0.5,1) -- (0.5,0);
			\draw[->,color=red] (0.25,1) -- (0.25,0);
			\draw[->] (0,0.5) -- (0.25,0.5);
			\draw[->] (1,0.5) -- (0.5,0.5);
			\end{tikzpicture} & &
			\begin{tikzpicture}[xscale=-0.5,yscale=0.5,rotate=180]
			\draw (0,0) rectangle (1,1);
			\draw[->] (0.5,1) -- (0.5,0);
			\draw[->,color=red] (0.25,1) -- (0.25,0);
			\draw[->] (0,0.5) -- (0.25,0.5);
			\draw[->] (1,0.5) -- (0.5,0.5);
			\end{tikzpicture}& \raisebox{1ex}{$\mapsto$} &
			\begin{tikzpicture}[xscale=-0.5,yscale=0.5,rotate=-270]
			\draw (0,0) rectangle (1,1);
			\draw[->] (0.5,1) -- (0.5,0);
			\draw[->,color=red] (0.25,1) -- (0.25,0);
			\draw[->] (0,0.5) -- (0.25,0.5);
			\draw[->] (1,0.5) -- (0.5,0.5);
			\end{tikzpicture} & &
			\begin{tikzpicture}[xscale=-0.5,yscale=0.5,rotate=-270]
			\draw (0,0) rectangle (1,1);
			\draw[->] (0.5,1) -- (0.5,0);
			\draw[->,color=red] (0.25,1) -- (0.25,0);
			\draw[->] (0,0.5) -- (0.25,0.5);
			\draw[->] (1,0.5) -- (0.5,0.5);
			\end{tikzpicture}& \raisebox{1ex}{$\mapsto$} &
			\begin{tikzpicture}[xscale=-0.5,yscale=0.5]
			\draw (0,0) rectangle (1,1);
			\draw[->] (0.5,1) -- (0.5,0);
			\draw[->,color=red] (0.25,1) -- (0.25,0);
			\draw[->] (0,0.5) -- (0.25,0.5);
			\draw[->] (1,0.5) -- (0.5,0.5);
			\end{tikzpicture} \\
			\begin{tikzpicture}[scale=0.5]
			\draw (0,0) rectangle (1,1);
			\draw[->] (0.5,1) -- (0.5,0);
			\draw[->,color=red] (0.25,1) -- (0.25,0);
			\draw[->] (0,0.5) -- (0.25,0.5);
			\draw[->] (1,0.5) -- (0.5,0.5);
			\end{tikzpicture}& \raisebox{1ex}{$\mapsto$} &
			\begin{tikzpicture}[scale=0.5,rotate=90]
			\draw (0,0) rectangle (1,1);
			\draw[->] (0.5,1) -- (0.5,0);
			\draw[->,color=red] (0.25,1) -- (0.25,0);
			\draw[->] (0,0.5) -- (0.25,0.5);
			\draw[->] (1,0.5) -- (0.5,0.5);
			\end{tikzpicture} & & 
			\begin{tikzpicture}[scale=0.5,rotate=90]
			\draw (0,0) rectangle (1,1);
			\draw[->] (0.5,1) -- (0.5,0);
			\draw[->,color=red] (0.25,1) -- (0.25,0);
			\draw[->] (0,0.5) -- (0.25,0.5);
			\draw[->] (1,0.5) -- (0.5,0.5);
			\end{tikzpicture}& \raisebox{1ex}{$\mapsto$} &
			\begin{tikzpicture}[scale=0.5,rotate=180]
			\draw (0,0) rectangle (1,1);
			\draw[->] (0.5,1) -- (0.5,0);
			\draw[->,color=red] (0.25,1) -- (0.25,0);
			\draw[->] (0,0.5) -- (0.25,0.5);
			\draw[->] (1,0.5) -- (0.5,0.5);
			\end{tikzpicture} & &
			\begin{tikzpicture}[scale=0.5,rotate=180]
			\draw (0,0) rectangle (1,1);
			\draw[->] (0.5,1) -- (0.5,0);
			\draw[->,color=red] (0.25,1) -- (0.25,0);
			\draw[->] (0,0.5) -- (0.25,0.5);
			\draw[->] (1,0.5) -- (0.5,0.5);
			\end{tikzpicture}& \raisebox{1ex}{$\mapsto$} &
			\begin{tikzpicture}[scale=0.5,rotate=270]
			\draw (0,0) rectangle (1,1);
			\draw[->] (0.5,1) -- (0.5,0);
			\draw[->,color=red] (0.25,1) -- (0.25,0);
			\draw[->] (0,0.5) -- (0.25,0.5);
			\draw[->] (1,0.5) -- (0.5,0.5);
			\end{tikzpicture} & &
			\begin{tikzpicture}[scale=0.5,rotate=270]
			\draw (0,0) rectangle (1,1);
			\draw[->] (0.5,1) -- (0.5,0);
			\draw[->,color=red] (0.25,1) -- (0.25,0);
			\draw[->] (0,0.5) -- (0.25,0.5);
			\draw[->] (1,0.5) -- (0.5,0.5);
			\end{tikzpicture}& \raisebox{1ex}{$\mapsto$} &
			\begin{tikzpicture}[scale=0.5]
			\draw (0,0) rectangle (1,1);
			\draw[->] (0.5,1) -- (0.5,0);
			\draw[->,color=red] (0.25,1) -- (0.25,0);
			\draw[->] (0,0.5) -- (0.25,0.5);
			\draw[->] (1,0.5) -- (0.5,0.5);
			\end{tikzpicture} \\
			
			\begin{tikzpicture}[scale=0.5]
			\draw (0,0) rectangle (1,1);
			\draw[<->] (0.5,1) -- (0.5,0);
			\draw[<->,color=red] (0.25,1) -- (0.25,0.25) -- (1,0.25);
			\draw[<->] (0,0.5) -- (1,0.5);
			\end{tikzpicture}& \raisebox{1ex}{$\mapsto$} &
			\begin{tikzpicture}[scale=0.5,rotate=90]
			\draw (0,0) rectangle (1,1);
			\draw[<->] (0.5,1) -- (0.5,0);
			\draw[<->,color=red] (0.25,1) -- (0.25,0.25) -- (1,0.25);
			\draw[<->] (0,0.5) -- (1,0.5);
			\end{tikzpicture} & & 
			\begin{tikzpicture}[scale=0.5,rotate=90]
			\draw (0,0) rectangle (1,1);
			\draw[<->] (0.5,1) -- (0.5,0);
			\draw[<->,color=red] (0.25,1) -- (0.25,0.25) -- (1,0.25);
			\draw[<->] (0,0.5) -- (1,0.5);
			\end{tikzpicture}& \raisebox{1ex}{$\mapsto$} &
			\begin{tikzpicture}[scale=0.5,rotate=180]
			\draw (0,0) rectangle (1,1);
			\draw[<->] (0.5,1) -- (0.5,0);
			\draw[<->,color=red] (0.25,1) -- (0.25,0.25) -- (1,0.25);
			\draw[<->] (0,0.5) -- (1,0.5);
			\end{tikzpicture} & &
			\begin{tikzpicture}[scale=0.5,rotate=180]
			\draw (0,0) rectangle (1,1);
			\draw[<->] (0.5,1) -- (0.5,0);
			\draw[<->,color=red] (0.25,1) -- (0.25,0.25) -- (1,0.25);
			\draw[<->] (0,0.5) -- (1,0.5);
			\end{tikzpicture}& \raisebox{1ex}{$\mapsto$} &
			\begin{tikzpicture}[scale=0.5,rotate=270]
			\draw (0,0) rectangle (1,1);
			\draw[<->] (0.5,1) -- (0.5,0);
			\draw[<->,color=red] (0.25,1) -- (0.25,0.25) -- (1,0.25);
			\draw[<->] (0,0.5) -- (1,0.5);
			\end{tikzpicture} & &
			\begin{tikzpicture}[scale=0.5,rotate=270]
			\draw (0,0) rectangle (1,1);
			\draw[<->] (0.5,1) -- (0.5,0);
			\draw[<->,color=red] (0.25,1) -- (0.25,0.25) -- (1,0.25);
			\draw[<->] (0,0.5) -- (1,0.5);
			\end{tikzpicture}& \raisebox{1ex}{$\mapsto$} &
			\begin{tikzpicture}[scale=0.5]
			\draw (0,0) rectangle (1,1);
			\draw[<->] (0.5,1) -- (0.5,0);
			\draw[<->,color=red] (0.25,1) -- (0.25,0.25) -- (1,0.25);
			\draw[<->] (0,0.5) -- (1,0.5);
			\end{tikzpicture} \\
			
			\begin{tikzpicture}[scale=0.5]
			\draw (0,0) rectangle (1,1);
			\draw[->] (0.5,1) -- (0.5,0);
			\draw[->] (0,0.5) -- (0.5,0.5);
			\draw[->] (1,0.5) -- (0.5,0.5);
			\draw[->,color=red] (0,0.25) -- (0.5,0.25);
			\draw[->,color=red] (1,0.25) -- (0.5,0.25);
			\end{tikzpicture}& \raisebox{1ex}{$\mapsto$} &
			\begin{tikzpicture}[scale=0.5,rotate=90]
			\draw (0,0) rectangle (1,1);
			\draw[->] (0.5,1) -- (0.5,0);
			\draw[->] (0,0.5) -- (0.5,0.5);
			\draw[->] (1,0.5) -- (0.5,0.5);
			\draw[->,color=red] (0,0.25) -- (0.5,0.25);
			\draw[->,color=red] (1,0.25) -- (0.5,0.25);
			\end{tikzpicture} & & 
			\begin{tikzpicture}[scale=0.5,rotate=90]
			\draw (0,0) rectangle (1,1);
			\draw[->] (0.5,1) -- (0.5,0);
			\draw[->] (0,0.5) -- (0.5,0.5);
			\draw[->] (1,0.5) -- (0.5,0.5);
			\draw[->,color=red] (0,0.25) -- (0.5,0.25);
			\draw[->,color=red] (1,0.25) -- (0.5,0.25);
			\end{tikzpicture}& \raisebox{1ex}{$\mapsto$} &
			\begin{tikzpicture}[scale=0.5,rotate=180]
			\draw (0,0) rectangle (1,1);
			\draw[->] (0.5,1) -- (0.5,0);
			\draw[->] (0,0.5) -- (0.5,0.5);
			\draw[->] (1,0.5) -- (0.5,0.5);
			\draw[->,color=red] (0,0.25) -- (0.5,0.25);
			\draw[->,color=red] (1,0.25) -- (0.5,0.25);
			\end{tikzpicture} & &
			\begin{tikzpicture}[scale=0.5,rotate=180]
			\draw (0,0) rectangle (1,1);
			\draw[->] (0.5,1) -- (0.5,0);
			\draw[->] (0,0.5) -- (0.5,0.5);
			\draw[->] (1,0.5) -- (0.5,0.5);
			\draw[->,color=red] (0,0.25) -- (0.5,0.25);
			\draw[->,color=red] (1,0.25) -- (0.5,0.25);
			\end{tikzpicture}& \raisebox{1ex}{$\mapsto$} &
			\begin{tikzpicture}[scale=0.5,rotate=270]
			\draw (0,0) rectangle (1,1);
			\draw[->] (0.5,1) -- (0.5,0);
			\draw[->] (0,0.5) -- (0.5,0.5);
			\draw[->] (1,0.5) -- (0.5,0.5);
			\draw[->,color=red] (0,0.25) -- (0.5,0.25);
			\draw[->,color=red] (1,0.25) -- (0.5,0.25);
			\end{tikzpicture} & &
			\begin{tikzpicture}[scale=0.5,rotate=270]
			\draw (0,0) rectangle (1,1);
			\draw[->] (0.5,1) -- (0.5,0);
			\draw[->] (0,0.5) -- (0.5,0.5);
			\draw[->] (1,0.5) -- (0.5,0.5);
			\draw[->,color=red] (0,0.25) -- (0.5,0.25);
			\draw[->,color=red] (1,0.25) -- (0.5,0.25);
			\end{tikzpicture}& \raisebox{1ex}{$\mapsto$} &
			\begin{tikzpicture}[scale=0.5]
			\draw (0,0) rectangle (1,1);
			\draw[->] (0.5,1) -- (0.5,0);
			\draw[->] (0,0.5) -- (0.5,0.5);
			\draw[->] (1,0.5) -- (0.5,0.5);
			\draw[->,color=red] (0,0.25) -- (0.5,0.25);
			\draw[->,color=red] (1,0.25) -- (0.5,0.25);
			\end{tikzpicture}		
			\end{array}\] 
			\caption{The relabeling map $\mathfrak{R}$ which replaces each tile with its corresponding rotation by $\frac{1}{2}\pi$.}
			\label{fig:rotation_tiles}
		\end{figure}
		
		It is easy to verify that $\rho$ and $\mu$ are valid extended symmetries, as they respect the conditions on the arrowheads and tails, and the sublattice comprised of only crosses. Also, we see that $\psi$ sends both $\rho$ and $\mu$ to the desired matrices, and that the mappings $\mathfrak{R}_\infty,\mathfrak{M}_\infty:\mathcal{A}^{\Z^2}\to\mathcal{A}^{\Z^2}$ commute with the corresponding rigid symmetries of the coordinate axes. Thus, $\langle\rho,\mu\rangle$ is a copy of $D_4$ contained in $\Sym(X_{\rm Rob},\Z^2)$, as desired.\qed
		
	\end{dem}
	
	We remark that the proof above used the structure of the Robinson shift $X_{\rm Rob}$ exclusively to compute the set of directions of fractures associated to this shift, and that extended symmetries preserve this set in other contexts as well. This suggests that this technique is open to generalization to other subshifts, even in higher dimensions, although possibly replacing the concept of ``direction of fracture'' with ``hyperplane of fracture'', as we need to separate half-spaces of $\Z^d$, whose boundaries are akin to $(d-1)$-dimensional affine spaces, albeit discrete. Thus, to propose a generalization of this concept we need a few definitions:
	
	\begin{defi}
		A \textbf{hyperplane} $H\subseteq\Z^d$ is a coset of a direct summand of $\Z^d$ of rank $d-1$; this is, $H$ is a nonempty subset of $\Z^d$ such that:
		\begin{enumerate}[label=(\arabic*)]
			\item $H=H_0+\vec{v}$ for some subgroup of $\Z^d$ with rank $d-1$ and some vector $\vec{v}\in\Z^d$, and
			\item there exists some other vector $\vec{w}\in\Z^d$ such that $\Z^d=H_0\oplus\Z\vec{w}$.
		\end{enumerate}
	\end{defi}
	Thus, we suggest the following tentative definition for a fracture in a $d$-di\-men\-sio\-nal subshift:
	\begin{defi}
		Let $X$ be a (strongly aperiodic) $\Z^d$-subshift. We say that $X$ has a \textbf{fracture} in the direction of the hyperplane $H=H_0+\vec{v}$ if for some $x\in X$ there are two half-spaces $S^+,S^-$ separated by $H$ (i.e. $S^+\cap S^-=S^+\cap H=S^-\cap H=\varnothing$) such that for some ``sufficiently large'' subset $B\subseteq H_0$ there is a family $\{x^{(\vec{b})}\}_{\vec{b}\in B}$ of points of $X$ such that:
		\[x^{(\vec{b})}|_{S^+}=x|_{S^+},\qquad x^{(\vec{b})}|_{S^-}=\sigma_{\vec{b}}(x)|_{S^-}.  \]
	\end{defi}
	Here, an appropiate definition of ``sufficiently large'' will depend on the subshift that is being studied. For instance, in the case of the Robinson shift we only needed $B$ to contain two points ($\{0,k_1\vec{q}\}$) for our argument due to the hierarchical structure of $X_{\rm Rob}$, albeit $B$ in this shift actually is an infinite set, $2\Z\vec{q}$. In all cases, as long as we apply a consistent restriction to the possible instances of $B$, we see that an extended symmetry $f$ must send a point of fracture to another point of fracture due to the generalized Curtis-Hedlund-Lyndon theorem, and thus $\psi(f)$ is a matrix that acts by permutation on the set of all hyperplanes of fracture of $X$. For sufficiently rigid shifts, this should result in a strong restriction on the matrix group $\psi[\Sym(X,\Z^d)]$.
	
	\section{The minimal case}
	
	The discussion above shows the key idea behind the showcased method: the hierarchical structure of the aforementioned subshifts forces the appearance of ``special directions'', which result in a geometrical invariant that needs to be preserved by extended symmetries. By identifying these special directions via combinatorial or dynamic properties, we can effectively restrict $\psi[\Sym(X,\Z^d)]$ enough to effectively compute it in terms of $\Aut(X,\Z^d)$.
	
	However, in the above discussion we focused specifically on the Robinson tiling $X_{\rm Rob}$ and the substitutive subshift $\shift{\theta}$ which (usually) are not minimal. To exhibit the above mentioned special directions, a key point was using certain points that exhibit ``fracture-like'' behavior, which are not present in the minimal subset of each of these subshifts. However, since the ``special directions'' come from the hierarchical structure of the subshift, they ought to be present in its minimal subset as well, and thus should impose the same restrictions on the set of extended symmetries. We proceed to show that this is actually the case, starting with the Robinson tiling as follows:
	
	\begin{cor}
		Let $M_{\rm Rob}\subset X_{\rm Rob}$ be the unique minimal subshift contained in $X_{\rm Rob}$. Then, $\Sym(M_{\rm Rob},\Z^2)\cong\Z^2\rtimes D_4$.
	\end{cor}
	
	\begin{dem}
		%\rotatebox{90}{\includegraphics[height=\heightof{M}]{tile1_mini.pdf}}
		Using the substitution rules devised by Gähler in \cite{Ga2013}, we can show that $M_{\rm Rob}$ contains a point, say $x$, that has only copies of the first tile from figure \ref{fig:robinson_tiles}, pointing to the right, on the horizontal strip $\Z\times\{0\}$, and corresponding tiles of the same kind pointing downwards in $\{0\}\times\Z^+$ and upwards in $\{0\}\times\Z^-$. Mirrored and rotated versions of this configuration exist as points of $M_{\rm Rob}$ as well (of which one specific rotation may be observed in Figure \ref{fig:robinson_supertile_separation}); a similar argument holds for the fifth tile.
		
		Any point from $\mathcal{H}=\overline{\{\sigma_{(n,0)}(x):n\in\Z \}}$ has the same horizontal strip of copies of the same tile on $\Z\times\{0\}$, and, due to the local rules of the Robinson tiling, any configuration with support $\Z\times[-n,n]$ from some point $y\in\mathcal{H}$ must be $(m,0)$-periodic for some sufficiently large $m$. Note that this $m$ must diverge to $\infty$ as $n\to\infty$, because no point from $M_{\rm Rob}$ has nontrivial periods.
		
		Let $f\in\Sym(M_{\rm Rob},\Z^2)$ be an extended symmetry. For any sufficiently large value of $k\in\N$, the window of this $f$ is contained in $\Z\times[-k,k]$. Thus, due to Theorem \ref{teo:curtishedlundlyndon}, we may choose a sufficiently large $k$ such that the image of $\Z\times[-k,k]$ under the matrix $\psi(f)$ contains the set $L_{a,b}(\widetilde{k})\dfn\{(u,v)\in\Z^2:-\tilde{k}\le au+bv\le\tilde{k} \}$ for any desired $\tilde{k}>0$ and some $a,b\in\Z$, and thus $y|_{\Z\times[-n,n]}$ determines $f(y)|_{L_{a,b}(\widetilde{k})}$ entirely. Since the strip $y|_{\Z\times[-k,k]}$ is periodic, the restriction $f(y)|_{L_{a,b}(\tilde{k})}$ must have a period as well, which we can choose as a multiple of $(-b,a)$.

		We may now either proceed with a combinatorial or dynamic argument; we show both as they are closely related, starting with the combinatorial method. For this, suppose that $ba\ne 0$, which implies that $\psi(f)$ maps $\vec{e}_1$ to a direction that is not parallel to the coordinate axes. Since the $n$-th order supertiles increase in size exponentially with $n$, and so do the associated ``square drawings'' determined by the crosses, the strip $L_{a,b}(\widetilde{k})$ must pass through the vertical lines (comprised of copies of rotations of the second, third, fourth or fifth tiles from Figure \ref{fig:robinson_tiles}) associated with the corresponding square of a $n$-th order supertile for all sufficiently large $n$ (as it is not parallel to any of the sides of such squares). Thus, this configuration cannot have a nontrivial period, since due to the positions of the $n$-th order supertiles such a period cannot have a horizontal or vertical component smaller than $2^n$, which applies for any sufficiently large $n$. We conclude, by this contradiction, that $\psi(f)$ maps $\vec{e}_1$ to a vector parallel to the coordinate axes; a similar argument holds with $\vec{e}_2$.

		From the dynamical perspective, we may proceed as in \cite{BRY2018} by employing the maximal equicontinuous factor (MEF) of the Robinson tiling. As shown in the aforementioned work by Gähler \cite{Ga2013}, the MEF of the Robinson tiling is a two-dimensional solenoid\footnote{The solenoid $\mathds{S}_p$ is the compact abelian group obtained as an inverse limit of the system $\R/\Z\gets\R/\Z\gets\R/\Z\gets\dotsc$, where each morphism is the mapping $x\mapsto px\pmod{1}$. A $d$-dimensional solenoid is defined analogously.} $\mathds{S}_2^2$, and the fiber of the corresponding mapping $\rho:M_{\rm Rob}\epi\mathds{S}_2^2$ is $28$-to-$1$ in the set of all points from $M_{\rm Rob}$ that are comprised of four infinite-order supertiles \cite{Ga2013}, including the above constructed $x$. By similar arguments to the ones from Corollary 1 of \cite{BRY2018}, an extended symmetry must map $x$, whose corresponding fiber $\rho^{-1}[\{\rho(x)\}]$ has $28$ different points, to another point with a corresponding fiber of cardinality $28$, comprised of four infinite-order supertiles. By employing the periodicity of the strip $x|_{\Z\times[-k,k]}$ and the local behavior of an extended symmetry, we can conclude that the matrix $\psi(f)$ maps $\Z\times[-k,k]$ to the support of the corresponding periodic strip in the image point, which (up to translation) must be either of the form $\Z\times[-\tilde{k},\tilde{k}]$ or $[-\tilde{k},\tilde{k}]\times\Z$; in both cases $\vec{e}_1$ must be mapped to a cardinal direction, as above.
		
		As the same holds for $\vec{e}_2$, we see that the matrix must be one of the eight matrices corresponding to the standard embedding of $D_4$ into $\GL_2(\Z)$, leading to the same conclusion as in the non-minimal case. The restrictions of the explicit mappings shown in the previous section to $M_{\rm Rob}$ constitute by themselves a copy of $D_4$ in $\Sym(M_{\rm Rob},\Z^2)$, as desired.\qed
	\end{dem}
	
	The analysis on the minimal subset of the Robinson tiling above suggests as well methods to study the minimal substitutive subshift obtained from a bijective substitution. We will prove the following result:
	
	\begin{teo}
		\label{teo:sym_grp_for_minimal_bij_subst}
		For a $d$-dimensional, nontrivial, primitive, bijective substitution $\theta$ over an alphabet $\Alf$ with faithful shift action, the following holds for the associated minimal substitutive subshift $\shift{\theta}^\circ$:
			\[\psi[\Sym(\shift{\theta}^\circ,\Z^d)]\le Q_d<\GL_d(\Z). \]
		Hence, the extended symmetry group is virtually $\Z^d$. 
	\end{teo}
	
	The previous result allows us to decompose $\Sym(\shift{\theta}^\circ,\Z^d)$ into a (semidirect) product of $\Z^d$ (the subgroup generated by the shifts), a subset of $S_{\abs{\Alf}}$ representing the relabeling maps, and a subset of $Q_d<\GL_d(\Z)$ corresponding to lattice transformations, in the same way as we did with $\Sym(\shift{\theta},\Z^d)$.
	
	We may proceed either combinatorially or dynamically, as above. For the combinatorial case, we need the following simple lemma:
	
	\begin{lem}
		\label{lem:substitutive_fracture}
		There exist two points $x^{(1)},x^{(2)}\in\shift{\theta}^\circ$ such that $x^{(1)}|_{\Z^+_0\times\Z^{d-1}}=x^{(2)}|_{\Z^+_0\times\Z^{d-1}}$, but $(x^{(1)})_{\vec{k}}\ne(x^{(2)})_{\vec{k}}$ for any $\vec{k}\in\Z^-\times\Z^{d-1}$.
	\end{lem}
	
	\begin{dem}
		Since the action $\Z^d\actson[\sigma]\shift{\theta}^\circ$ is faithful and minimal, there must be symbols $a,b,c\in\Alf$, with $b\ne c$, such that, for some points $x,y\in\shift{\theta}^\circ$, $x_{\vec{0}}=y_{\vec{0}}=a$ but $x_{-\vec{e}_1}=b,y_{-\vec{e}_1}=c$. If this were not the case, for any point $x\in\shift{\theta}^\circ$ the symbol $x_{\vec{k}}$ would determine $x_{\vec{k}+\vec{e}_1}$ uniquely; since $\abs{\Alf}<\infty$ this would result in a direction of periodicity shared by all points in $\shift{\theta}^\circ$, a contradiction.
		
		As usual, we may replace $\theta$ by $\theta^m$ for a sufficiently large $m$ such that every periodic point of $\theta$ is a fixed point of $\theta^m$. By the previous observation, there exist two fixed points $x',y'\in\shift{\theta}^\circ$ such that $x'_{\vec{0}}=y'_{\vec{0}}=a$ and $x'_{-\vec{e}_1}=b,y'_{-\vec{e}_1}=c$. Since $x'$ and $y'$ are fixed points of the substitution, these symbols determine the corresponding quadrants entirely, and thus $x'$ and $y'$ match on the subset $Q_{\vec{\ind}}=(\Z^+_0)^d$ but (due to bijectiveness) differ in every symbol from $\Z^-\times(\Z^+_0)^{d-1}$. Taking $\vec{k}=(0,1,1,\dotsc,1)$, any ordered pair $(x^{(1)},x^{(2)})$ that is a limit point of the sequence $(\sigma_{\vec{k}}^m(x'),\sigma_{\vec{k}}^m(y'))_{m\ge 0}$ (note that such a pair exists by compactness) satisfies the desired condition. \qed
	\end{dem}
	
	Lemma \ref{lem:substitutive_fracture} provides us an analogue to the fracture points from the Robinson tiling; namely, it provides a separating hyperplane $\{0\}\times\Z^{d-1}$ that ``splits'' $\Z^d$ into two half-spaces $S^-$ and $S^+$, and two points $x,y$ which match on one half-space, say $S^+$, but not the other. An analogue of this result applies to any other hyperplane of the form $\Z^{k-1}\times\{0\}\times\Z^{d-k},1\le k\le d$ or any of their translates; thus, the same argument used in the case of the Robinson tiling $X_{\rm Rob}$ applies here to show that the set of all ``fracture hyperplanes'' has to be preserved. Hence, Theorem \ref{teo:sym_grp_for_minimal_bij_subst} follows immediately from the following result:
	
	\begin{lem}
		\label{lem:fracture_norm_dir}
		For the minimal subshift $\shift{\theta}^\circ$ given by a bijective substitution $\theta$ under the above hypotheses, call $\vec{v}\in\Z^d\setminus\{\vec{0}\}$ a \textbf{fracture normal direction} if there is some $N>0$ and two disjoint subsets $S^+,S^-\sqsubset\Z^d$ of the form:
			\[S^\pm\dfn\{\vec{k}\in\Z^d: {\pm\langle\vec{k},\vec{v}\rangle}\ge N \}, \]
		such that, for some $x,y\in\shift{\theta}^\circ$, $x|_{S^+}=y|_{S^+}$ but $x_{\vec{k}}\ne y_{\vec{k}}$ for any $\vec{k}\in S^-$. Then, the set of fracture normal directions of $\shift{\theta}^\circ$ is $\{h\vec{e}_j: h\in\Z\setminus\{0\},1\le j\le d \}$.
	\end{lem}
	
	In other words, for $\shift{\theta}^\circ$ as in Lemma \ref{lem:fracture_norm_dir} the set of all possible fracture hyperplanes is exactly the set of translations of coordinate hyperplanes of the form $\Z^{k-1}\times\{0\}\times\Z^{d-k}$. The proof is similar to the argument above for the non-minimal case:
	
	\begin{dem}
		As stated above, Lemma \ref{lem:substitutive_fracture} shows that the set of fracture normal directions of $\shift{\theta}^\circ$ contains $\{h\vec{e}_j: h\in\Z\setminus\{0\},1\le j\le d \}$. Suppose $\vec{v}$ is an additional fracture normal direction not contained in this set; as it is not parallel to the coordinate axes, an argument similar to the one from Lemma \ref{lem:finite_lines} shows that the set:
			\[L_N \dfn \{\vec{k}\in\Z^d: \abs{\langle\vec{k},\vec{v}\rangle}<N \}=\Z^d\setminus (S^+\dotcup S^-), \]
		has finite intersection with some $\Z\vec{e}_j,1\le j\le d$, and the size of this intersection is bounded by a value depending linearly on $N$ and the entries of $\vec{v}$. By Lemma \ref{lem:codified_system_subst}, any pair of fracture points $x,y$ associated to the direction $\vec{v}$ can be written as the concatenation of patterns of the form $\theta^m(a),a\in\Alf$, whose supports are rectangles with side length depending exponentially on $m$. Thus, by choosing a sufficiently large $m$, the support $R$ of one of these patterns $\theta^m(a)$ has nonempty intersection with both $S^+$ and $S^-$. 
		
		Since $x|_{S^+}=y|_{S^+}$ and the substitution is bijective, we must have $x|_R=y|_R$. However, this implies $x|_{R\cap S^{-}}=y|_{R\cap S^{-}}$, where $R\cap S^{-}$ is a nonempty set, contradicting our hypothesis (as $x_{\vec{k}}\ne y_{\vec{k}}$ for all $\vec{k}\in S^{-}$). Thus, this $\vec{v}$ cannot be a fracture normal direction.\qed
	\end{dem}
	
	%Alternatively, we could proceed dynamically by emulating the procedure from \cite{BRY2018} (as mentioned above in the case of the Robinson tiling); as stated above in Lemma \ref{lem:integer_images_in_MEF}, if $\varphi:\shift{\theta}\epi\Z_{\vec{s}}$ is the factor map from the substitutive subshift (with the hypothesis above) to its maximal equicontinuous factor (which is a product of odometers), a point $x\in X$ such that $\varphi(x)\in\Z^d$ must be in the orbit of a $\theta$-periodic point, and, in particular, $\varphi(x)=\vec{0}$ if and only if $x$ is a $\theta$-periodic point.
	
	The same arguments as above allow us to conclude that $\psi(f)$ necessarily is a matrix from the standard copy of $Q_d$ in $\GL_d(\Z)$, for any $f\in\Sym(\shift{\theta}^\circ,\Z^d)$. Alternatively, we could argue dynamically in the same vein as \cite{BRY2018}, as follows. Let $\varphi:\shift{\theta}^\circ\epi\Z_{\vec{n}}$ be the mapping from the minimal substitutive subshift $\shift{\theta}^\circ$ to its maximal equicontinuous factor. As the number of $\theta$-periodic points is exactly the same as the number of admissible patterns from $\shift{\theta}$ with support $\{-1,0\}^d$, say $\ell=\abs{\lang_{\{0,1\}^d}(\shift{\theta}^\circ)}$, the mapping $\varphi$ is $\ell$-to-$1$ in the set of all $\theta$-periodic points $\Per_\theta(\shift{\theta}^\circ)$. An argument similar to the one employed in the proof of Lemma \ref{lem:integer_images_in_MEF} can be used to show the following result:
	
	\begin{lem}
		Under the above hypotheses, given $x\in\shift{\theta}$, define $J=\{j\in\{1,\dotsc,d\}:\varphi(x)_j\in\Z \}$. Then, $\Z^d$ can be partitioned into sets of the form:
			\[S_H\dfn\{\vec{n}\in\Z^d: n_j\ge\varphi(x)_j \text{ if }j\in H \land n_j<\varphi(x)_j\text{ if }j\in J\setminus H \}, \]
		for all $H\subseteq J$, such that any point $y\in\shift{\theta}$ with $\varphi(x)=\varphi(y)$ is entirely determined by the configuration $y|_U$, where $U\subset\Z^d$ is any subset with non-trivial intersection with all of the $S_H$.
	\end{lem}
	
	In particular, we may choose $U$ as a translation of the rectangle $R_J$ with support $\{-1,0\}^J\times\{0\}^{\{1,\dots,d\}\setminus J}$. This imposes a strong restriction on the points $y$ with the same image as $x$ (namely, that there can be at most $\abs{\lang_{R_J}(\shift{\theta})}$ such points) and, by faithfulness, this shows that $\varphi$ cannot be $\ell$-to-$1$ in $\shift{\theta}\setminus\Per_\theta(\shift{\theta})$. By the argument exhibited in $\cite{BRY2018}$, extended symmetries preserve the cardinality of the fibers in the MEF and thus must map $\Per_\theta(\shift{\theta})$ to itself; the desired result then follows from the local behavior of an extended symmetry, via a similar analysis as in the case of the Robinson tiling above.
	
	\section*{Acknowledgements}
	
	The author would like to thank his thesis advisor, PhD Michael H. Schraudner, for his guidance, advice, insight and ideas with regard to the development of the current work, and Michael Baake and Franz Gähler from Bielefeld University for their helpful comments and advice regarding further developments and additions. The author would also like to kindly acknowledge the financial support received from CONICYT Doctoral Fellowship 21171061, year 2017.
	
	The final section of this paper was written during a stay at Bielefeld, which was financed by the CRC 1283 at Bielefeld University. Additional financial support for this stay was received from University of Chile.
	
	%The final section of this work was worked upon during a visit to the University of Bielefeld, which was jointly financed by the University of Chile and project SFB 1283
	
	\bibliographystyle{hplain}
	
	\bibliography{coven_highdim}
\end{document}